\documentclass{amsart}

\newtheorem{theorem}{Theorem}[section]
\newtheorem{lemma}[theorem]{Lemma}
\newtheorem{proposition}[theorem]{Proposition}
\newtheorem{corollary}[theorem]{Corollary}

\theoremstyle{definition}
\newtheorem{definition}[theorem]{Definition}
\newtheorem{example}[theorem]{Example}
\newtheorem{remark}[theorem]{Remark}

\usepackage{amscd,amssymb}

\begin{document}

\title[Normal and Normally Outer Automorphisms]
{Normal and Normally Outer Automorphisms\\
of Free Metabelian Nilpotent Lie algebras }

\author[\c{S}ehmus F\i nd\i k]
{\c{S}ehmus F\i nd\i k}

\address{Department of Mathematics,
\c{C}ukurova University, 01330 Balcal\i,
 Adana, Turkey}
\email{sfindik@cu.edu.tr}

\thanks
{The research of the author was partially supported by the
 Council of Higher Education (Y\"OK) in Turkey}

\subjclass[2010]
{17B01, 17B30, 17B40.}
\keywords{free Lie algebras, free metabelian nilpotent Lie algebras, inner automorphisms, normal automorphisms.}

\begin{abstract} 
Let $L_{m,c}$ be the free $m$-generated metabelian nilpotent of class $c$ Lie algebra
over a field of characteristic 0. An automorphism $\varphi$ of $L_{m,c}$ is called normal if
$\varphi(I)=I$ for every ideal $I$ of the algebra $L_{m,c}$.
Such automorphisms form a normal subgroup $\text{\rm N}(L_{m,c})$
of $\text{\rm Aut}(L_{m,c})$ containing the group of inner automorphisms.
We describe the group of normal automorphisms of $L_{m,c}$
and the quotient group of $\text{\rm Aut}(L_{m,c})$ modulo $\text{\rm N}(L_{m,c})$.
\end{abstract}

\maketitle

\section*{Introduction}
Let $L_m$ be the free $m$-generated Lie algebra over a field $K$ of characteristic 0, $m\geq 2$, and let $L_{m,c}=L_m/(L_m''+L_m^{c+1})$ be the free $m$-generated metabelian nilpotent of class $c$
Lie algebra. This is the relatively free algebra of rank $m$ in the variety of Lie
algebras $\mathfrak A^2\cap \mathfrak N_c$, where $\mathfrak A^2$ is the
metabelian (solvable of class 2) variety of Lie algebras and ${\mathfrak N}_c$ is the variety of
all nilpotent Lie algebras of class at most $c$.

An automorphism of an algebra is called normal if it preserves every ideal of the algebra.
Similarly, an automorphism of a group is normal if it preserves every normal subgroup of the group.
Such automorphisms form a normal subgroup of the group of all automorphisms.
The goal of our paper is to describe the group of normal automorphisms $\text{\rm N}(L_{m,c})$
and the quotient group $\text{\rm Aut}(L_{m,c})/\text{\rm N}(L_{m,c})$ of normally outer automorphisms of
the Lie algebra $L_{m,c}$.
The corresponding problem for the group of normal automorphisms of free metabelian nilpotent groups
was studied by Endimioni \cite{E1,E2,E3}. He showed that the normal automorphisms $\theta$ of
a free metabelian nilpotent group $G$ are exactly the atomorphisms of the form
\[
\theta(x)=x(x,u_1)^{k(1)}\ldots(x,u_m)^{k(m)},
\]
where $u_1,\ldots , u_m$ are elements of $G$, the exponents $k(1),\ldots, k(m)$ are integers. (As
usual, the commutator $(a, b)$ in the group case is defined by $(a,b)= a^{-1}b^{-1}ab$.).
Endimioni also proved that the group of normal automorphisms of
free metabelian nilpotent group $G$ is metabelian, generalizing a result of Gupta \cite{G}
for the group of $\text{\rm IA}$-automorphisms in a two-generated metabelian group. Initially, automorphisms
of the form $\theta(x)=x(x,u_1)^{k(1)}\ldots(x,u_m)^{k(m)}$ were studied by Kuzmin \cite{K}.

The group of normal automorphisms of free groups has been studied by Lubotzky \cite{L}.
He showed that $\text{\rm N}(G)=\text{\rm Inn}(G)$, for any finitely generated free group $G$.
Lue \cite{Lu} gave a short proof of this fact using the Freiheitssatz for groups established by Magnus \cite{M}.
The Freiheitssatz for Lie algebras was proved by Shirshov \cite{S}. Makar-Limanov \cite{Ma} proved it
for associative algebras over a field of characteristic zero. 
Following the idea of Lue \cite{Lu} we show that the free Lie algebra $L_m$ does not have nontrivial
normal automorphisms for any $m\geq 2$ and over a field of any characteristic.
For the proof we apply the Freiheitssatz for Lie algebras and use the Hopf property of free Lie algebras.
The same result holds for free associative algebras over a field of characteristic $0$.
The key step of the proof was suggested to us by Ualbai Umirbev.
If we replace $L_m$ with a relatively free algebra in a proper subvariety of all Lie algebras
it may happen that many normal automorphisms appear.
In particular, this holds for the free metabelian nilpotent
Lie algebra $L_{m,c}$. Since every inner automorphism of $L_{m,c}$ is normal, the algebra
$L_{m,c}$ posseses nontrivial normal automorphisms.                                                                                                                  
These automorphisms form a normal subgroup of $\text{\rm Aut}(L_{m,c})$.

Our first main result is similar to the result of Endimioni \cite{E1,E3} in the case of groups
but there are some essential differences.
We show that the group of normal automorphisms is included in the subgroup $\text{\rm IA}(L_{m,c})$ of the automorphisms
which induce the identity map modulo the commutator ideal
of $L_{m,c}$ when $m\geq3,c\geq2$ or $m=2,c\geq4$. 
In the exceptional cases, i.e. $(m,c)=(2,2)$ or $(m,c)=(2,3)$, every normal automorphism
acts on the generators of $L_{m,c}$ as a nonzero scalar times an $\text{\rm IA}$-automorphism.
For the proof we define a special type of automorphisms called generalized inner automorphisms and we 
describe the group of normal automorphisms in terms of them.
We also show that the group of normal automorphisms $\text{\rm N}(L_{m,c})$
is an abelian group when $m\geq3,c=2$,
is a nilpotent of class 2 group when $m\geq3,c=3$ and
is a metabelian group when $m\geq2,c\geq4$ or $(m,c)=(2,2)$. Finally, $\text{\rm N}(L_{m,c})$
is a nilpotent of class two--by--abelian group when $(m,c)=(2,3)$
which is an analogue of the result of Gupta \cite{G} and Endimioni \cite{E2}.

A result of Shmel'kin \cite{Sh} states that the free metabelian Lie algebra $F_m=L_m/L''_m$ can be embedded into the abelian
wreath product $A_m\text{\rm wr}B_m$, where $A_m$ and $B_m$ are $m$-dimensional abelian Lie algebras
with bases  $\{a_1,\ldots,a_m\}$ and $\{b_1,\ldots,b_m\}$, respectively.
The elements of $A_m\text{\rm wr}B_m$ are of the form 
\[
\sum_{i=1}^ma_if_i(t_1,\ldots,t_m)+\sum_{i=1}^m\beta_ib_i,
\]
where the $f_i$'s are polynomials in $K[t_1,\ldots,t_m]$ and $\beta_i\in K$. This allows to introduce partial derivatives
in $F_m$ with values in $K[t_1,\ldots,t_m]$ and the Jacobian matrix $J(\phi)$ of an endomorphism $\phi$ of $F_m$.
Restricted on the semigroup $\text{\rm IE}(F_m)$ of endomorphisms of $F_m$ which are identical modulo
the commutator ideal $F_m'$, the map $J:\phi\to J(\phi)$ is a semigroup monomorphism of $\text{\rm IE}(F_m)$
into the multiplicative semigroup of the algebra $M_m(K[t_1,\ldots,t_m])$ of $m\times m$ matrices
with entries from $K[t_1,\ldots,t_m]$.
In the present work we consider the embedding of the free metabelian nilpotent Lie algebra $L_{m,c}$ into
the wreath product $A_m\text{\rm wr}B_m$ modulo the ideal $(A_m\text{\rm wr}B_m)^{c+1}$. The automorphism
group $\text{\rm Aut}(L_{m,c})$ is a semidirect
product of the normal subgroup $\text{\rm IA}(L_{m,c})$ and the general linear group $\text{\rm GL}_m(K)$.
Considering the group $\text{\rm IN}(L_{m,c})$ of normal $\text{\rm IA}$-automorphisms,
for the description of the factor group
$\Gamma\text{\rm N}(L_{m,c})=\text{\rm Aut}(L_{m,c})/\text{\rm N}(L_{m,c})$ it is sufficient to know only
$\text{\rm IA}(L_{m,c})/\text{\rm IN}(L_{m,c})$.
Drensky and F\i nd\i k \cite{DF} gave the explicit form of the Jacobian matrices of
the coset representatives of the outer automorphisms in
$\text{\rm IA}(L_{m,c})/\text{\rm Inn}(L_{m,c})$. Since $\text{\rm Inn}(L_{m,c})$ is included in
the group of normal automorphisms, $\text{\rm IA}(L_{m,c})/\text{\rm IN}(L_{m,c})$
is a homomorphic image of $\text{\rm IA}(L_{m,c})/\text{\rm Inn}(L_{m,c})$
and we find explicitely coset representatives of $\text{\rm IN}(L_{m,c})$.

The paper is organized as follows.
In the first section, we introduce normal and normally outer automorphisms and discuss the
relations between $\text{\rm N}(L_{m,c})$ and the normal subgroup $\text{\rm IA}(L_{m,c})$. We also discuss
the normal automorphisms of the free Lie algebra $L_m$.
In the second section we define the group of
generalized inner automorphisms and give necessary information
about its group structure. In the third section we describe the group of normal
automorphisms in terms of the group of generalized inner automorphisms. Finally we give the explicit
form of the Jacobian matrices of the normal automorphisms and of the Jacobian matrices of the coset representatives of
normally outer $\text{\rm IA}$-automorphisms. We also give the explicit form of the Jacobian matrices
of the coset representatives of the normal automorphisms modulo the group of inner automorphisms $\text{\rm Inn}(L_{m,c})$.

\section{Preliminaries}
Let $L_m$ be the free Lie algebra of rank $m\geq 2$ over a field $K$ of characteristic 0
with free generators $y_1,\ldots,y_m$ and let $L_{m,c}=L_m/(L_m''+L_m^{c+1})$ be the free
metabelian nilpotent of class $c$ Lie algebra freely generated by
$x_1,\ldots,x_m$, where $x_i=y_i+(L_m''+L_m^{c+1})$, $i=1,\ldots,m$.
We use the commutator notation for the Lie multiplication. Our commutators are left normed:
\[
[u_1,\ldots,u_{n-1},u_n]=[[u_1,\ldots,u_{n-1}],u_n],\quad n=3,4,\ldots.
\]
In particular,
\[
L_{m,c}^k=\underbrace{[L_{m,c},\ldots,L_{m,c}]}_{k\text{ times}}.
\]
For each $v\in L_{m,c}$, the linear operator $\text{\rm ad}v:L_{m,c}\to L_{m,c}$ defined by
\[
u(\text{\rm ad}v)=[u,v],\quad u\in L_{m,c},
\]
is a derivation of $L_{m,c}$ which is nilpotent and $\text{\rm ad}^cv=0$
because $L_{m,c}^{c+1}=0$.
Hence the linear operator
\[
\exp(\text{\rm ad}v)=1+\frac{\text{\rm ad}v}{1!}+\frac{\text{\rm ad}^2v}{2!}+\cdots
=1+\frac{\text{\rm ad}v}{1!}+\frac{\text{\rm ad}^2v}{2!}+\cdots+\frac{\text{\rm ad}^{c-1}v}{(c-1)!}
\]
is well defined and it is an inner automorphism  of $L_{m,c}$.
The set of all such automorphisms form a normal subgroup $\text{\rm Inn}(L_{m,c})$
of the group of all automorphisms $\text{\rm Aut}(L_{m,c})$ of $L_{m,c}$.

Let $\varphi$ be an automorphism of an algebra $R$ such that 
$\varphi(I)=I$ for every ideal $I$ of the algebra $R$. Such automorphisms 
are called $normal$ automorphisms which we denote by $\text{\rm N}(R)$. Clearly these 
automorphisms form a normal subgroup of the group of all automorphisms $\text{\rm Aut}(R)$ of $R$.
The factor group $\text{\rm Aut}(R)/\text{\rm N}(R)$
is the group of $normally$ $outer$ (or $N$-outer) automorphisms and is denoted by $\Gamma\text{\rm N}(R)$.

The next lemma gives the form of normal automorphisms of $L_{m,c}$.

\begin{lemma}\label{normalLmc}
Let $\varphi$ be a normal automorphism of the algebra $L_{m,c}$. Then $\varphi$ is of the form
\[
\varphi: x_i\to \alpha x_i+\sum_{j=1}^{m}[x_i,x_j]f_{ij}(\text{\rm ad}x_1,\ldots,\text{\rm ad}x_m),
\quad i=1,\ldots,m,\quad \alpha\in K^*,
\]
where $f_{ij}(t_1,\ldots,t_m)\in K[t_1,\ldots,t_m]$ and $K^*$ is the set of nonzero elements
of the field $K$.
\end{lemma}

\begin{proof}
Let $\varphi$ be a normal automorphism of the algebra $L_{m,c}$. Hence $\varphi$ induces
a normal automorphism $\overline\varphi$ of the abelian algebra $\overline L_{m,c}=L_{m,c}/L'_{m,c}$.
The automorphism group of $\overline L_{m,c}$ coincides with the general linear group $GL_m(K)$ and
the normal automorphisms of $\overline L_{m,c}$ are the elements of $GL_m(K)$ which preserve
the vector subspaces of $\overline L_{m,c}$. Applying to the vector subspace $K\overline x_i$ we obtain that
$\overline\varphi(\overline x_i)=\alpha_i\overline x_i$, $\alpha_i\in K^*$.
Similarly, for $i\neq j$,
\begin{align}
\overline\varphi(\overline {x_i+x_j})&=\alpha_i\overline x_i+\alpha_j\overline x_j\nonumber\\
&=\beta(\overline x_i+\overline x_j),\quad \beta\in K^*.\nonumber
\end{align}
Thus $\alpha_i=\alpha_j=\alpha$. 
Hence $\varphi$ has the form
\[
\varphi: x_i\to \alpha x_i+u_i, \quad \alpha\in K^*,\, u_i\in L'_{m,c},\, i=1,\ldots,m.
\]

It is well known in a metabelian Lie algebra $G$, see e.g. \cite{Ba}, that
\[
[v_1,v_2,v_{\sigma(3)},\ldots,v_{\sigma(k)}]
=[v_1,v_2,v_3,\ldots,v_k],\quad v_1,\ldots,v_k\in G,
\]
where $\sigma$ is an arbitrary permutation of $3,\ldots,k$, i.e. the operators $\text{\rm ad}v$ , $v\in G$ commute
when acting on $G'$. The vector space $L_{m,c}'$ has a basis consisting of all
\[
[x_{i_1},x_{i_2},x_{i_3},\ldots,x_{i_k}],\quad 1\leq i_j\leq m,\, i_1>i_2\leq i_3\leq\cdots\leq i_k,\, k\leq c,
\]
and we may permute the elements $x_{i_3},\ldots,x_{i_k}$. Reordering the elements $x_1,\ldots,x_m$ by
\[
x_i<x_1<\cdots<x_{i-1}<x_{i+1}<\cdots<x_m
\]
we obtain that the subspace of $L_{m,c}'$ spanned by the commutators essentially depending on $x_i$, has a basis
\[
[x_i,x_j,x_{i_3},\ldots,x_{i_k}],\quad j\neq i,\, 1\leq i_3\leq\cdots\leq i_k, \,k\leq c.
\]
Hence the normal automorphism $\varphi$ of $L_{m,c}$ has the form
\[
\varphi: x_i\to \alpha x_i+\sum_{j\neq i}[x_i,x_j]f_{ij}(\text{\rm ad}x_1,\ldots,\text{\rm ad}x_m)+g_i(\hat x_i), 
\]
where $\alpha\in K^*$, $f_{ij}(t_1,\ldots,t_m)\in K[t_1,\ldots,t_m]$ and $g_i(\hat x_i)\in L'_{m,c}$ does not depend
on $x_i$.

For a fixed $i=1,\ldots,m$ let us consider the ideal $J_i$ of $L_{m,c}$ generated by the element $x_i$, 
Since $\varphi$ is normal and $\varphi(x_i)\in J_i$ we obtain that
\[
g_i(\hat x_i)\in J_i,\quad i=1,\ldots,m,
\]
and hence
\[
g_i(\hat x_i)=0,\quad i=1,\ldots,m,
\]
because every element in $J_i$ depends on $x_i$. Thus we have
\[
\varphi(x_i)=\alpha x_i+\sum_{j\neq i}[x_i,x_j]f_{ij}(\text{\rm ad}x_1,\ldots,\text{\rm ad}x_m),
\]
which completes the proof.
\end{proof}

A similar proof holds also for the free Lie algebra $L_m$. But we use
the fact that, applying the anticommutativity and the Jacobian identity,
linear combinations of commutators of $L_m$ depending essentially on $y_i$
can be rewritten as linear combinations of left normed commutators of the form
\[
[y_i,y_{i_2},\ldots,y_{i_k}],\quad y_{i_2},\ldots,y_{i_k}\in\{y_1,\ldots,y_m\}.
\]

\begin{lemma}\label{freelie}
Let $\varphi$ be a normal automorphism of the algebra $L_m$. Then $\varphi$ is of the form
\[
\varphi: y_i\to \alpha y_i+y_if_i(\text{\rm ad}Y), \quad i=1,\ldots,m,\quad \alpha\in K^*,
\]
where $f_i(\text{\rm ad}Y)=f_i(\text{\rm ad}y_1,\ldots,\text{\rm ad}y_m)$ and every polyomial
$f_i(t_1,\ldots,t_m)$, $i=1,\ldots,m$, belongs to the free associative algebra $K\langle t_1,\ldots,t_m\rangle$.
\end{lemma}

Recall that an algebra $R$ is $Hopfian$, if it cannot be mapped onto itself with nontrivial kernel.
The following fact is folklorely known.

\begin{lemma}\label{hopfLie}
Finitely generated free Lie algebras and free associative algebras over any field of
arbitrary characteristic are Hopfian.
\end{lemma}

For example it is stated for relatively free algebras of finite rank as Exercise 4.10.21,
page 137 in the book of Bahturin \cite{Ba}.  The proof is similar to the proof in the group case,
see Section 4.1 of the book by Neumann \cite{N}, and repeats the steps of the proof of Theorem 9, page 104 \cite{Ba}.
The proof of \cite [Exercise 4.10.21]{Ba} uses only the fact that over an infinite field relatively
free algebras $F(\mathfrak U)$ are graded and that $\cap_{m\geq1}F^m(\mathfrak U)=0$
which is obviously true for free Lie algebras and free associative algebras over any field.

The analogue of the Freiheitssatz in group theory \cite{Ma}
was proved by Shirshov \cite{S} in the case of Lie algebras in any characteristic.
For associative algebras it was obtained by Makar-Limanov \cite{Ma}
when characteristic of the base field is $0$. The problem is still open for associative
algebras over a field of positive characteristic (see e.g.
the book by Bokut' and Kukin \cite{BK}). We shall state the result for free Lie algebras only.

\begin{theorem}\label{freedomalgebra}{\rm(Shirshov \cite{S})}
Let $L(Y)$ be the free Lie algebra freely generated
by $Y=\{ y_1,\ldots,y_m\}$.
If $f(Y)\in L(Y)$ does not belong to the subalgebra
generated by $ y_1,\ldots,y_{m-1}$, then $(f(Y))\cap L(y_1,\ldots,y_{m-1})=0$ where
$(f(Y))$ is the ideal of $L(Y)$ generated by $f(Y)$.
\end{theorem}

The idea to use the Freiheitssatz in the following proof was suggested to us by Ualbai Umirbaev.

\begin{corollary}\label{coordinate}
If every monomial of $f(Y)\in L(Y)$ depends on $y_m$ and
 $f(Y)\notin L(y_m)=Ky_m$, then
$f(Y)$ is not an image of $y_m$ under an automorphism of the algebra
$L(Y)$, i.e. $f(Y)$ is not a coordinate.
\end{corollary}

\begin{proof}
Let $\varphi$ be an automorphism of $L(Y)$ and let  $\varphi(y_m)=f(Y)$, i.e. $f=f(Y)$ be a coordinate.
Clearly $f\in (y_m)\vartriangleleft L(Y)$ because every monomial of $f$ depends on $y_m$.
Let $y_m\notin (f)\vartriangleleft L(Y)$. 
This means that $f$ depends also on the variables $y_1,\ldots,y_{m-1}$.
Since $\varphi: L(Y)\to L(Y)$ is an automorphism and $\varphi(y_m)=f$, then
\[
L(Y)/(f)\cong L(\varphi(y_1),\ldots,\varphi(y_{m-1}))\cong L(y_1,\ldots,y_{m-1}).
\]
On the other hand $L(Y)/(y_m)\cong L(y_1,\ldots,y_{m-1})$.
As a result
\[
L(Y)/(f)\cong L(Y)/(y_m)\cong L(y_1,\ldots,y_{m-1}).
\]
Let us consider the natural homomorphism
\[
\pi:  L(y_1,\ldots,y_{m-1})\cong L(Y)/(f)\to L(Y)/(y_m)\cong L(y_1,\ldots,y_{m-1}).
\]
$\pi$ is onto and $\ker\pi\neq0$ because $y_m\notin (f)$. But $L(y_1,\ldots,y_{m-1})$
has the $Hopf$ $property$ by Lemma \ref{hopfLie}. This is in contradiction with $\ker\pi\neq0$.
Then $(f)=(y_m)\vartriangleleft L(Y)$. But $f$ depends also on the other variables,
for example, without loss of generality, depends also on $y_1$. Applying Theorem \ref{freedomalgebra} (the Freiheitssatz)
we get that $y_2,\ldots,y_m$ generate a free algebra of rank $m-1$ in $L(Y)/(f)$.
But this is not true for $L(Y)/(y_m)$, because $\overline y_m\neq \overline 0$ in $L(Y)/(f)$ while  $\overline y_m= \overline 0$ in $L(Y)/(y_m)$.
\end{proof}

Lubotzky \cite{L} showed that the group of normal automorphisms of a free group $G$
is equal to the group of inner automorphisms of $G$, i.e. $\text{\rm N}(G)=\text{\rm Inn}(G)$
and Lue \cite{Lu} gave an alternative proof of the statement.
Our next theorem states that the free Lie algebra $L_m$ does not have nontrivial
normal automorphisms for any $m\geq 2$.
The idea of the proof is similar to the idea of the proof of the paper by Lue \cite{Lu} for free groups. 

\begin{theorem}\label{normaloffreeLie}
Let $L_m$ be the free Lie algebra of rank $m\geq2$ over a field $K$ of characteristic $0$
with free generators $y_1,\ldots,y_m$. Then $L_m$ does not have nontrivial normal automorphisms.
\end{theorem}

\begin{proof}
Applying Lemma \ref{freelie} and Corollary \ref{coordinate}
we see that $L_m$ does not have nontrivial normal $\text{\rm IA}$-automorphisms,
i.e. the normal automorphisms of $L_m$ are of the form $y_i\to \alpha y_i$, $i=1,\ldots,m$, $\alpha\in K^*$. 

Let $m=2$ and consider the ideal $I$ of $L_2$ generated by $f=y_1-[y_1,y_2]$.
Let $\varphi$ be a normal automorphism of $L_2$ of the form
\[
y_1\to \alpha y_1,\quad y_2\to \alpha y_2,\quad \alpha\in K^*,
\]
and assume that $\alpha\neq1$. Since $\varphi$ is normal
\[
\varphi(f)=\alpha y_1-\alpha^2[y_1,y_2]\in I.
\]
Hence we have the system
\begin{align}
y_1-[y_1,y_2]&\equiv0 \quad (\text{\rm mod }I)\nonumber\\
\alpha y_1-\alpha^2[y_1,y_2]&\equiv0 \quad (\text{\rm mod }I)\nonumber
\end{align}
Since $\alpha\neq0,1$, then $y_1\equiv0 \quad (\text{\rm mod }I)$
and $[y_1,y_2]\equiv0\quad(\text{\rm mod}I)$ which means that $I=(y_1)\vartriangleleft L_2$.
Now consider the ideal $J$ of $L_2$ generated by all commutators
$u\in L_2$ such that $\text{\rm deg}_{y_1}(u)\geq2$.
Then
\[
\overline L_2=L_2/J=\text{\rm span}\{\overline y_2,[\overline y_1,\underbrace{\overline y_2,\ldots,\overline y_2}_{k}]\mid k\geq0\}.
\]
Recall that $f=y_1-[y_1,y_2]$. Clearly $[\overline f,\overline y_1]=\overline0$ in $\overline L_2$.
So the elements of $\overline I$ in $\overline L_2$ are linear combinations of
\[
u_k=[\overline f,\underbrace{\overline y_2,\ldots,\overline y_2}_{k}], \quad k\geq0,.\]
Thus 
\[
\overline I=\left\{ \sum_{k\geq0}\beta_ku_k\mid\beta_k\in K,\,\sum_{k\geq0}\beta_k=0\right\}.
\]
This means that
$\overline y_1\in(\overline y_1)$ while $\overline y_1\notin\overline I$, because the only coefficient of $\overline y_1$ is $1$.
Thus $(\overline y_1)\neq\overline I$
which is in contradiction with $I=(y_1)$. Hence $\alpha=1$.

Now let $m\geq 3$ and let $\varphi$ be a normal automorphism of $L_m$ of the form 
\[
\varphi(y_i)=\alpha y_i\quad i=1,\ldots,m,\quad \alpha\in K^*.
\]
We consider the ideal $I$ of $L_m$ generated by the elements $y_1-[y_1,y_2],y_3\ldots,y_m$. 
Since $\varphi$ is normal then $\varphi(I)=I$ and it induces a normal automorphism of
$L_m/I$ which is isomorphic to $L_2/(y_1-[y_1,y_2])$.
But we already know that in this case $\alpha=1$.
\end{proof}

\begin{remark}
An analogue of Theorem \ref{normaloffreeLie} holds for
free associative algebras $K\langle Y\rangle= K\langle y_1,\ldots,y_m\rangle$ over a field of characteristic 0.
Repeating the main steps of the proof of Theorem \ref{normaloffreeLie} we obtain that the only possibility is that
$f(Y)$ depends on $y_m$ only. We extend $\varphi$ to an
automorphism $\overline\varphi$ of the algebra $\overline K\langle y_1,\ldots,y_m\rangle$ where $\overline K$
is the algebraic closure of $K$. If $\text{\rm deg}(f(Y)=f(y_m))=d>1$, then
\[
\varphi(y_m)=f(y_m)=a_0(y_m-\alpha_1)\cdots(y_m-\alpha_d),\quad a_0,\alpha_1,\ldots,\alpha_d\in\overline K,
\]
is a product of several polynomials which is impossible:

Applying $\varphi^{-1}$ we obtain that the degree of 
\[
y_m=a_0(\varphi^{-1}(y_m)-\alpha_1)\cdots(\varphi^{-1}(y_m)-\alpha_d)
\]
is bigger than 1.
\end{remark}

The situation in the case of free metabelian nilpotent
Lie algebra $L_{m,c}$ is different.
Applying Lemma \ref{normalLmc} it is easy to see that
$\text{\rm Inn}(L_{m,c})\subset \text{\rm N}(L_{m,c})$.
Hence $L_{m,c}$ posseses nontrivial normal automorphisms. The group $\text{\rm N}(L_{m,c})$ is not necessarily included in
the normal subgroup $\text{\rm IA}(L_{m,c})$ of $\text{\rm Aut}(L_{m,c})$ of the automorphisms which induce the identity
map modulo the commutator ideal of $L_{m,c}$. Our next lemma states that in some cases
$\text{\rm N}(L_{m,c})\subset\text{\rm IA}(L_{m,c})$.

\begin{lemma}\label{nia}
$\text{\rm (i)}$ If $m\geq3$ and $c=2$, then $\text{\rm N}(L_{m,2})\subset \text{\rm IA}(L_{m,2})$.

$\text{\rm (ii)}$ If $m\geq3$ and $c=3$, then $\text{\rm N}(L_{m,3})\subset \text{\rm IA}(L_{m,3})$.

$\text{\rm (iii)}$ If $m\geq2$ and $c\geq4$, then $\text{\rm N}(L_{m,c})\subset \text{\rm IA}(L_{m,c})$.
\end{lemma}

\begin{proof}$\text{\rm (i)}$ Let $\varphi$ be a normal automorphism of $L_{m,2}$, $m\geq3$. By Lemma \ref{normalLmc}
$\varphi$ has the form
\[
\varphi: x_i\to \alpha x_i+\sum_{j=1}^{m}\beta_{ij}[x_i,x_j], \quad i=1,\ldots,m,\quad \alpha\in K^*, 
\]
where $\beta_{ij}\in K$.
Let us consider the ideal $J$ generated by $u=x_1+[x_2,x_3]$. $J$ has a basis 
\[
x_1+[x_2,x_3],\, [x_1,x_j],\, j=2,\ldots,m.
\]
Since $\varphi$ is normal
\[
\varphi(u)=\alpha x_1+\alpha^2[x_2,x_3]+\beta_{12}[x_1,x_2]+\cdots+\beta_{1m}[x_1,x_m]\in J
\]
Clearly the summand $\alpha x_1+\alpha^2[x_2,x_3]$ is included in the vector space spanned by the element
$u$. Thus $\alpha=\alpha^2$ or $\alpha=1$.

$\text{\rm (ii)}$ Let $\varphi$ be a normal automorphism of $L_{m,3}$, $m\geq3$. By Lemma \ref{normalLmc}
$\varphi$ has the form
\[
\varphi: x_i\to \alpha x_i+\sum_{j=1}^{m}[x_i,x_j]f_{ij}, \quad i=1,\ldots,m,\quad \alpha\in K^*, 
\]
where $f_{ij}\in K[\text{\rm ad}x_1,\ldots,\text{\rm ad}x_m]$.
We can express $\varphi$ as
\[
\varphi: x_i\to \alpha x_i+\sum_{j=1}^{m}[x_i,x_j](f_{ij,0}+f_{ij,1}),
\]
where $f_{ij,0}\in K$, $f_{ij,1}\in \omega/\omega^2$. Here $\omega$ states for
the augmentation ideal of $K[\text{\rm ad}x_1,\ldots,\text{\rm ad}x_m]$.

Let us consider the ideal $J$ generated by $u=[x_1,x_2]+[x_1,x_3,x_3]$. $J$ has a basis 
\[
[x_1,x_2]+[x_1,x_3,x_3],\, [x_1,x_2,x_j],\, j=1,\ldots,m.
\]
Since $\varphi$ is normal, $\varphi(u)\in J$. Easy calculations give that
\[
\varphi(u)=\alpha^2[x_1,x_2]+\alpha^3[x_1,x_3,x_3]-\alpha\sum_{j=1}^{m}[x_2,x_j,x_1]f_{2j,0}
+\alpha\sum_{j=1}^{m}[x_1,x_j,x_2]f_{1j,0}.
\]
Clearly the summand $\alpha^2[x_1,x_2]+\alpha^3[x_1,x_3,x_3]$ is included in the vector space spanned by the element
$[x_1,x_2]+[x_1,x_3,x_3]$. Thus $\alpha^2=\alpha^3$ or $\alpha=1$.

$\text{\rm (iii)}$ Let $\varphi$ be a normal automorphism of $L_{m,c}$, $m\geq2,\,c\geq4$. By Lemma \ref{normalLmc}
$\varphi$ has the form
\[
\varphi: x_i\to \alpha x_i+\sum_{j=1}^{m}[x_i,x_j]f_{ij}, \quad i=1,\ldots,m,\quad \alpha\in K,
\]
where $f_{ij}\in K[\text{\rm ad}x_1,\ldots,\text{\rm ad}x_m]$.

Let us consider the ideal $J$ generated by 
$v=\underbrace{[x_1,x_2,\ldots,x_2]}_{c-1}+\underbrace{[x_1,x_2,x_1,\ldots,x_1]}_{c}$.
$J$ has a basis consisting of $v$ and the elements of the form
\[
\underbrace{[x_1,x_2,\ldots,x_2,x_j]}_{c},\quad j=1,\ldots,m.
\]
Since $\varphi$ is normal $\varphi(v)\in J$. Similar steps as $\text{\rm (ii)}$ give that
\[
\alpha^{c-1}\underbrace{[x_1,x_2,\ldots,x_2]}_{c-1}+\alpha^c\underbrace{[x_1,x_2,x_1,\ldots,x_1]}_{c}
\]
is included in 
the vector space spanned by the element $v$. Thus $\alpha=1$.
\end{proof}

Let $F_m=L_m/L''_m$ be the free metabelian Lie algebra
of rank $m$. We shall denote the free generators of $F_m$ with the same symbols $x_1,\ldots,x_m$ as the free generators of
$L_{m,c}$, but now $x_i=y_i+L''_m$, $i=1,\ldots,m$.
Let $K[t_1,\ldots,t_m]$ be the
(commutative) polynomial algebra over $K$ freely generated by the
variables $t_1,\ldots,t_m$ and let $A_m$ and $B_m$ be the abelian
Lie algebras with bases $\{a_1,\ldots,a_m\}$ and
$\{b_1,\ldots,b_m\}$, respectively. Let $C_m$ be the free right
$K[t_1,\ldots,t_m]$-module with free generators $a_1,\ldots,a_m$.
We give it the structure of a Lie algebra with trivial multiplication.
The abelian wreath product
$A_m\text{\rm wr}B_m$ is equal to the semidirect sum $C_m\leftthreetimes B_m$. The elements
of $A_m\text{\rm wr}B_m$ are of the form
$\sum_{i=1}^ma_if_i(t_1,\ldots,t_m)+\sum_{i=1}^m\beta_ib_i$, where
$f_i$ are polynomials in $K[t_1,\ldots,t_m]$ and $\beta_i\in K$.
The multiplication in $A_m\text{\rm wr}B_m$ is defined by
\[
[C_m,C_m]=[B_m,B_m]=0,
\]
\[
[a_if_i(t_1,\ldots,t_m),b_j]=a_if_i(t_1,\ldots,t_m)t_j,\quad i,j=1,\ldots,m.
\]
Hence $A_m\text{\rm wr}B_m$ is a metabelian Lie algebra and every mapping $\{x_1,\ldots,x_m\}\to A_m\text{\rm wr}B_m$
can be extended to a homomorphism $F_m\to A_m\text{\rm wr}B_m$. In particular, as a special case of the embedding theorem of Shmel'kin \cite{Sh},
the mapping $x_i\to a_i+b_i$, $i=1,\ldots,m$, can be extended
to an embedding of $F_m$ into $A_m\text{\rm wr}B_m$.

Both $F_m$ and $A_m\text{\rm wr}B_m$ are graded algebras. The monomials in $A_m\text{\rm wr}B_m$
of degree 1 are of the form $a_i,b_j$ and of degree $n\geq2$ have the form 
$a_it_1\ldots t_{n-1}$. Let us consider the ideal $(A_m\text{\rm wr}B_m)^{c+1}$ spanned by
the elements of $A_m\text{\rm wr}B_m$ of length at least $c+1$. Then the quotient $(A_m\text{\rm wr}B_m)/(A_m\text{\rm wr}B_m)^{c+1}$ is metabelian and nilpotent of class $c$
and the homomorphism $\varepsilon: L_{m,c}\to (A_m\text{\rm wr}B_m)/(A_m\text{\rm wr}B_m)^{c+1}$ defined by
$\varepsilon(x_i)=a_i+b_i$, $i=1,\ldots,m$, is a monomorphism. If
\[
f=\sum[x_i,x_j]f_{ij}(\text{\rm ad}x_1,\ldots,\text{\rm ad}x_m),\quad f_{ij}(t_1,\ldots,t_m)\in K[t_1,\ldots,t_m]/\Omega^c,
\]
where $\Omega$ is
the augmentation ideal of $K[t_1,\ldots,t_m]$, then
\[
\varepsilon(f)=\sum(a_it_j-a_jt_i)f_{ij}(t_1,\ldots,t_m).
\]
The next lemma follows from \cite{Sh}, see also \cite {BD}.

\begin{lemma}\label{metabelian rule}
The element $\sum_{i=1}^ma_if_i(t_1,\ldots,t_m)$ of
$C_m$ belongs to $\varepsilon ( L_{m,c}')$ if and only if
$\sum_{i=1}^mt_if_i(t_1,\ldots,t_m)\equiv0\quad (\text{\rm mod }\Omega^{c+1})$.
\end{lemma}

The embedding of $L_{m,c}$ into
$A_m\text{\rm wr}B_m/(A_m\text{\rm wr}B_m)^{c+1}$ allows to introduce partial derivatives
in $L_{m,c}$ with values in $K[t_1,\ldots,t_m]/\Omega^c$. If $f\in L_{m,c}$ and
\[
\varepsilon(f)=\sum_{i=1}^m\beta_ib_i+\sum_{i=1}^ma_if_i(t_1,\ldots,t_m),\quad \beta_i\in K,f_i\in K[t_1,\ldots,t_m]/\Omega^c,
\]
then
\[
\frac{\partial f}{\partial x_i}=f_i(t_1,\ldots,t_m).
\]
The Jacobian matrix $J(\phi)$ of an endomorphism $\phi$ of $L_{m,c}$
is defined as
\[
J(\phi)=\left(\frac {\partial \phi({x_j})}{\partial x_i}\right)
=\left(\begin{matrix}
\frac {\partial\phi({x_1})}{\partial x_1}&\cdots&\frac {\partial \phi({x_m})}{\partial x_1}\\
\vdots&\ddots&\vdots\\
\frac {\partial\phi({x_1})}{\partial x_m}&\cdots&\frac {\partial \phi({x_m})}{\partial x_m}\\
\end{matrix}\right)\in M_m(K[t_1,\ldots,t_m]/\Omega^c),
\]
where $M_m(K[t_1,\ldots,t_m]/\Omega^c)$ is the associative algebra of $m\times m$ matrices with entries from
$K[t_1,\ldots,t_m]/\Omega^c$. Let $\text{\rm IE}(L_{m,c})$ be the multiplicative semigroup of all endomorphisms
of $L_{m,c}$ which are identical modulo the commutator ideal $L_{m,c}'$.
Let $I_m$ be the identity $m\times m$ matrix and let $S$ be the subspace of
$M_m(K[t_1,\ldots,t_m]/\Omega^c)$ defined by
\[
S=\left \{(f_{ij})\in M_m(K[t_1,\ldots,t_m]/\Omega^c) \mid 
\sum_{i=1}^mt_if_{ij}\equiv 0(\text{\rm mod}\Omega^{c+1}),j=1,\ldots,m\right \}.
\]
Clearly $I_m+S$ is a subsemigroup of the multiplicative group of $M_m(K[t_1,\ldots,t_m]/\Omega^c)$.
If $\phi\in \text{\rm IE}(L_{m,c})$, then $J(\phi)=I_m+(s_{ij})$,
where $s_{ij}\in S$.
It is easy to check that if $\phi,\psi \in \text{\rm IE}(L_{m,c})$ then $J(\phi \psi)=J(\phi)J(\psi)$.
The following proposition is well known, see e.g. \cite {BD}.

\begin{proposition}\label{met}
The map $J:\text{\rm IE}(L_{m,c})\to I_m+S$ defined by
$\phi\to J(\phi)$ is an isomorphism of the semigroups $\text{\rm IE}(L_{m,c})=\text{\rm IA}(L_{m,c})$ and $I_m+S$.
\end{proposition}

Now we know that the group of normal automorphisms
$\text{\rm N}(L_{m,c})$ is included in the subgroup
$\text{\rm IA}(L_{m,c})$  when $m\geq3,c=3$ or $m\geq2,c\geq4$ and in other cases every normal automorphism
is a nonzero scalar times an $\text{\rm IA}$-automorphism.
The automorphism group $\text{\rm Aut}(L_{m,c})$ is a semidirect
product of the normal subgroup $\text{\rm IA}(L_{m,c})$ and the general linear group $\text{\rm GL}_m(K)$.
Considering the group of normal $\text{\rm IA}$-automorphisms $\text{\rm IN}(L_{m,c})$,
for the description of the factor group
$\Gamma\text{\rm N}(L_{m,c})=\text{\rm Aut}(L_{m,c})/\text{\rm N}(L_{m,c})$ it is sufficient to know only
$\text{\rm IA}(L_{m,c})/\text{\rm IN}(L_{m,c})$. Drensky and F\i nd\i k \cite{DF} gave the explicit form of the Jacobian matrices of
the coset representatives of the outer automorphisms in
$\text{\rm IA}(L_{m,c})/\text{\rm Inn}(L_{m,c})$. 
Since $\text{\rm Inn}(L_{m,c})$ is included in
the group of normal automorphisms, $\text{\rm IA}(L_{m,c})/\text{\rm IN}(L_{m,c})$ is the homomorphic image of the factor group $\text{\rm IA}(L_{m,c})/\text{\rm Inn}(L_{m,c})$.

\begin{lemma}\label{formula of out}{\rm(Drensky and F\i nd\i k \cite{DF})}
The automorphisms with the following Jacobian matrices are
coset representatives of the subgroup $\text{\rm Inn}(L_{m,c})$
of the group $\text{\rm IA}(L_{m,c})$:
\[
J(\theta)=I_m+\left(\begin{array}{llll}
s(t_2,\ldots,t_m)&f_{12}&\cdots&f_{1m}\\
t_1q_2(t_2,t_3,\ldots,t_m)+r_2(t_2,\ldots,t_m)&f_{22}&\cdots&f_{2m}\\
t_1q_3(t_3,\ldots,t_m)+r_3(t_2,\ldots,t_m)&f_{32}&\cdots&f_{3m}\\
\ \ \ \ \ \ \ \vdots&\ \ \vdots&\ \ddots&\ \ \vdots\\
t_1q_m(t_m)+r_m(t_2,\ldots,t_m)&f_{m2}&\cdots&f_{mm}\\
\end{array}\right),
\]
where $s,q_i,r_i,f_{ij}\in \Omega/\Omega^c$, i.e., are polynomials of degree $\leq c-1$ without constant terms.
They satisfy the conditions
\[
s+\sum_{i=2}^mt_iq_i\equiv 0,\quad \sum_{i=2}^mt_ir_i\equiv 0,\quad \sum_{i=1}^mt_if_{ij}\equiv 0
\quad (\text{\rm mod }\Omega^{c+1}), \quad j=2,\ldots,m,
\]
$r_i=r_i(t_2,\ldots,t_m)$, $i=1,\ldots,m$, does not depend on $t_1$, $q_i(t_i,\ldots,t_m)$,
$i=2,\ldots,m$, does not depend on $t_1,\ldots,t_{i-1}$
and  $f_{12}$ does not contain a summand $dt_2$, $d\in K$.
\end{lemma}

\section{Generalized Inner Automorphisms}

In this section we introduce a special type of automorphisms of the free metabelian nilpotent Lie algebra $L_{m,c}$.
We shall use these automorphisms in order to describe the group of normal automorphisms $\text{\rm N}(L_{m,c})$
of $L_{m,c}$ in the next section.

\begin{definition}
An automorphism $\psi$ of the algebra $L_{m,c}$  is called $generalized$ $inner$ $automorphism$ if $\psi$ has the form
\[
\psi: x_i\to x_i+\sum_{j=1}^{m}[x_i,x_j]f_j, \quad i=1,\ldots,m,
\]
where $f_j\in K[\text{\rm ad}x_1,\ldots,\text{\rm ad}x_m]$.
\end{definition}

One can check that every inner automorphism is a generalized inner automorphism. We give necessary
information for the structure of generalized inner automorphisms in the next lemmas and theorems.

\begin{lemma}\label{multiplicationginner} Let $\psi$ and $\phi$ be generalized inner automorphisms of the form
\[
\psi: x_i\to x_i+\sum_{j=1}^{m}[x_i,x_j]f_j, \quad i=1,\ldots,m,
\]
\[
\phi: x_i\to x_i+\sum_{t=1}^{m}[x_i,x_t]g_t, \quad i=1,\ldots,m,
\]
where $f_j,g_t\in K[\text{\rm ad}x_1,\ldots,\text{\rm ad}x_m]$. Then
the composition $\psi\phi$ is of the form
\[
\psi\phi: x_i\to x_i+\sum_{t=1}^{m}[x_i,x_t]g_t+\sum_{j=1}^{m}[x_i,x_j]f_j+
\sum_{j,t=1}^{m}[x_i,x_t,x_j]g_tf_j, \quad i=1,\ldots,m.
\]
\end{lemma}

\begin{proof} Let $\psi$ and $\phi$ be as above. Then
\begin{align}
\psi(\phi(x_i))&=\psi(x_i)+\sum_{t=1}^{m}[\psi(x_i),\psi(x_t)]g_t\nonumber\\
&=\psi(x_i)+\sum_{t=1}^{m}\left[x_i+\sum_{j=1}^{m}[x_i,x_j]f_j,x_t+\sum_{j=1}^{m}[x_t,x_j]f_j\right]g_t\nonumber\\
&=\psi(x_i)+\sum_{t=1}^{m}[x_i,x_t]g_t+\sum_{t,j=1}^{m}([x_i,x_j,x_t]-[x_t,x_j,x_i])f_jg_t\nonumber\\
&=x_i+\sum_{t=1}^{m}[x_i,x_j]g_t+\sum_{j=1}^{m}[x_i,x_j]f_j+
\sum_{j,t=1}^{m}[x_i,x_t,x_j]g_tf_j, \quad i=1,\ldots,m\nonumber.
\end{align}
\end{proof}

\begin{theorem}\label{ginnergroup}
Generalized inner automorphisms form a subgroup of the automorphism group $\text{\rm Aut}(L_{m,c})$.
\end{theorem}

\begin{proof}
Let $\psi$ and $\phi$ be generalized inner automorphisms of the form
\begin{align}
\psi: x_i\to x_i+\sum_{j=1}^{m}[x_i,x_j]f_j, \quad i=1,\ldots,m,\nonumber\\
\phi: x_i\to x_i+\sum_{t=1}^{m}[x_i,x_t]g_t, \quad i=1,\ldots,m\nonumber,
\end{align}
where $f_j,g_t\in K[\text{\rm ad}x_1,\ldots,\text{\rm ad}x_m]$.
Applying Lemma \ref{multiplicationginner} we have that
 
\begin{align}
(\psi\phi)(x_i)&=x_i+\sum_{t=1}^{m}[x_i,x_t]g_t+\sum_{j=1}^{m}[x_i,x_j]f_j+
\sum_{j,t=1}^{m}[x_i,x_t,x_j]g_tf_j\nonumber\\
&=x_i+\sum_{j=1}^{m}[x_i,x_j](g_j+f_j)+\sum_{t=1}^{m}[x_i,x_t]g_t\sum_{j=1}^{m}\text{\rm ad}x_jf_j\nonumber,
\end{align}
for every $i=1,\ldots,m$. Let us put $h_t=g_t\sum_{j=1}^{m}\text{\rm ad}x_jf_j$, $t=1,\ldots,m$. So we have
\begin{align}
(\psi\phi)(x_i)&=x_i+\sum_{j=1}^{m}[x_i,x_j](g_j+f_j)+\sum_{t=1}^{m}[x_i,x_t]h_t\nonumber\\
&=x_i+\sum_{j=1}^{m}[x_i,x_j](g_j+f_j+h_j)\nonumber\\
&=x_i+\sum_{j=1}^{m}[x_i,x_j]F_j\nonumber,
\end{align}
where $F_j=g_j+f_j+h_j$, $j=1,\ldots,m$. Thus the composition $\psi\phi$ is a generalized inner automorphism.
It remains to prove that for any inverse automorphism $\psi^{-1}$ of a generalized inner automorphism $\psi$
is also a generalized inner automorphism.
For this purpose it suffices to construct for each
integer $n\geq2$ a generalized inner automorphism $\psi_n$ such that $\psi_n\psi$ is of the form
\[
\psi\psi_n: x_i\to x_i+\sum_{j=1}^{m}[x_i,x_j]h_j, \quad i=1,\ldots,m,
\]
where $h_j \in \omega^{n-1}$, where $\omega$ states for
the augmentation ideal of $K[\text{\rm ad}x_1,\ldots,\text{\rm ad}x_m]$,
i.e. the length of the commutator $[x_i,x_j]h_j$ is at least $n+1$.
Let $\psi$ be of the form
\[
\psi: x_i\to x_i+\sum_{j=1}^{m}[x_i,x_j](f_{j0}+\cdots+f_{j,c-2}), 
\]
where $f_{j0}\in K$, $f_{jk}\in\omega^k/\omega^{k+1}$, $k=1,\ldots,c-2$.
Let us consider the generalized inner automorphism
\[
\psi_2:x_i\to x_i-\sum_{j=1}^{m}f_{j0}[x_i,x_j],\quad f_{j0}\in K.
\]
From Lemma \ref{multiplicationginner} we obtain that
\[
\psi\psi_2: x_i\to x_i+\sum_{j=1}^{m}[x_i,x_j](g_{j1}+\cdots+g_{j,c-2}), \quad g_{jk}\in\omega^k/\omega^{k+1}.
\]
Now consider the generalized inner automorphism 
\[
\psi_3:x_i\to x_i-\sum_{j=1}^{m}[x_i,x_j]g_{j1},\quad g_{j1}\in\omega.
\]
Similarly we have that
\[
\psi\psi_2\psi_3: x_i\to x_i+\sum_{j=1}^{m}[x_i,x_j](h_{j2}+\cdots+h_{j,c-2}), \quad h_{jk}\in\omega^k/\omega^{k+1}.
\]
Repeating this process we construct $\psi_2,\psi_3,\ldots,\psi_c$ and obtain that
\[
\psi\psi_2\psi_3\ldots\psi_c=1.
\]
\end{proof}

\begin{lemma}\label{psi(u)}
Let $\psi$ be a generalized inner automorphism of the form
\[
\psi: x_i\to x_i+\sum_{j=1}^{m}[x_i,x_j]f_j, \quad i=1,\ldots,m,
\]
where $f_j\in K[\text{\rm ad}x_1,\ldots,\text{\rm ad}x_m]$. Then for every $u\in L_{m,c}$
\[
\psi(u)=u+\sum_{j=1}^{m}[u,x_j]f_j.
\]
\end{lemma}

\begin{proof}
By linearity it is sufficient to show for every $k=1,2,\ldots$ that
\[
\psi([x_{i1},\ldots,x_{ik}])=[x_{i1},\ldots,x_{ik}]+\sum_{j=1}^{m}[[x_{i1},\ldots,x_{ik}],x_j]f_j.
\]
We make induction on the degree 
$k$ of the commutators. The case $k=1$ is trivial.
It is true for $k=2$:
\begin{align}
\psi[x_p,x_q]&=[\psi(x_p),\psi(x_q)]\nonumber\\
&=\left[x_p+\sum_{j=1}^{m}[x_p,x_j]f_j,x_q+\sum_{j=1}^{m}[x_q,x_j]f_j\right]\nonumber\\
&=[x_p,x_q]+\sum_{j=1}^{m}([x_p,x_j,x_q]-[x_q,x_j,x_p])f_j\nonumber\\
&=[x_p,x_q]+\sum_{j=1}^{m}[[x_p,x_q],x_j]f_j\nonumber.
\end{align}
Now assume that the equality holds for $k-1$. Then
\begin{align}
\psi([x_{i_1},\ldots,x_{i_k}])&=[\psi([x_{i_1},\ldots,x_{i_{k-1}}]),\psi(x_{i_k})]\nonumber\\
&=\left[[x_{i_1},\ldots,x_{i_{k-1}}]+\sum_{j=1}^{m}[[x_{i_1},
\ldots,x_{i_{k-1}}],x_j]f_j,x_{i_k}+\sum_{j=1}^{m}[x_{i_k},x_j]f_j\right]\nonumber\\
&=[x_{i_1},\ldots,x_{i_k}]+\sum_{j=1}^{m}[[x_{i_1},\ldots,x_{i_k}],x_j]f_j\nonumber.
\end{align}
\end{proof}

\begin{corollary}
The group of generalized inner automorphisms $\text{\rm GInn}(L_{m,c})$ is a subgroup of
the group of normal automorphisms $\text{\rm N}(L_{m,c})$.
\end{corollary}

\begin{proof} Let $\psi$ be a generalized inner automorphism of the form
\[
\psi: x_i\to x_i+\sum_{j=1}^{m}[x_i,x_j]f_j, \quad i=1,\ldots,m,
\]
where $f_j\in K[\text{\rm ad}x_1,\ldots,\text{\rm ad}x_m]$. Let $u$ be an element of an ideal 
$J$ of the free metabelian nilpotent Lie algebra $L_{m,c}$. From Lemma \ref{psi(u)} we know that
\[
\psi(u)=u+\sum_{j=1}^{m}[u,x_j]f_j.
\]
Hence $\psi(u)\in J$.
\end{proof}

Now we describe the group structure of the group of generalized inner automorphisms $\text{\rm GInn}(L_{m,c})$.

\begin{theorem}\label{gigroup}
$\text{\rm (i)}$ The group $\text{\rm GInn}(L_{m,2})$ is abelian;\\
$\text{\rm (ii)}$ The group $\text{\rm GInn}(L_{m,3})$ is nilpotent of class 2;\\
$\text{\rm (iii)}$ The group $\text{\rm GInn}(L_{m,c})$, $c\geq4$, is metabelian.
\end{theorem}

\begin{proof}
$\text{\rm (i)}$ Let $\psi,\phi\in\text{\rm GInn}(L_{m,2})$ be generalized inner automorphisms of the form
\[
\psi: x_i\to x_i+\sum_{j=1}^{m}\alpha_j[x_i,x_j], \quad i=1,\ldots,m,
\]
\[
\phi: x_i\to x_i+\sum_{j=1}^{m}\beta_j[x_i,x_j], \quad i=1,\ldots,m,
\]
where $\alpha_j,\beta_j\in K$ for $j=1,\ldots,m$. Then the composition $\psi\phi$ is
\begin{align}
\psi(\phi(x_i))&=\psi(x_i+\sum_{j=1}^{m}\beta_j[x_i,x_j])\nonumber\\
&=x_i+\sum_{j=1}^{m}\alpha_j[x_i,x_j]+\sum_{t=1}^{m}\beta_j[x_i,x_j])\nonumber\\
&=x_i+\sum_{j=1}^{m}(\alpha_j+\beta_j)[x_i,x_j].\nonumber
\end{align}
Thus $\psi\phi=\phi\psi$.

$\text{\rm (ii)}$ Let $\varphi,\phi,\gamma\in\text{\rm GInn}(L_{m,3})$
be generalized inner automorphisms of the form
\[
\varphi: x_i\to x_i+\sum_{j=1}^{m}[x_i,x_j]f_j, \quad i=1,\ldots,m,
\]
\[
\phi: x_i\to x_i+\sum_{j=1}^{m}[x_i,x_j]g_j, \quad i=1,\ldots,m,
\]
\[
\gamma: x_i\to x_i+\sum_{j=1}^{m}[x_i,x_j]h_j, \quad i=1,\ldots,m,
\]
where $f_j,g_j,h_j\in K[\text{\rm ad}x_1,\ldots,\text{\rm ad}x_m]$ and let
\[
u=\sum_{j=1}^{m}\text{\rm ad}x_jf_j,\quad v=\sum_{j=1}^{m}\text{\rm ad}x_jg_j,
\quad w=\sum_{j=1}^{m}\text{\rm ad}x_jh_j.
\]
Using the arguments of Theorem \ref{ginnergroup} we have that
\[
\varphi^{-1}=\varphi_2\varphi_3;\quad\phi^{-1}=\phi_2\phi_3;\quad\gamma^{-1}=\gamma_2\gamma_3,
\]
where
\[
\varphi_2: x_i\to x_i-\sum_{j=1}^{m}[x_i,x_j]f_j;\quad \varphi_3: x_i\to x_i+\sum_{j=1}^{m}[x_i,x_j]f_ju;
\]
\[
\phi_2: x_i\to x_i-\sum_{j=1}^{m}[x_i,x_j]g_j;\quad \phi_3: x_i\to x_i+\sum_{j=1}^{m}[x_i,x_j]g_jv;
\]
\[
\gamma_2: x_i\to x_i-\sum_{j=1}^{m}[x_i,x_j]h_j;\quad \gamma_3: x_i\to x_i+\sum_{j=1}^{m}[x_i,x_j]h_jw.
\]
Using Lemma \ref{multiplicationginner} direct calculations give that
\[
(\varphi,\phi)=\varphi^{-1}\phi^{-1}\varphi\phi=\varphi_2\varphi_3\phi_2\phi_3\varphi\phi
\]
has the form
\[
(\varphi,\phi):x_i\to x_i+\sum_{j=1}^{m}[x_i,x_j](g_ju-f_jv).
\]
Let us define $t=\sum_{j=1}^{m}\text{\rm ad}x_j(g_ju-f_jv)$. Then we obtain that
\[
(\varphi,\phi,\gamma):x_i\to x_i+\sum_{j=1}^{m}[x_i,x_j](h_jt-(g_ju-f_jv)w).
\]
Since the polynomials $h_jt,(g_ju-f_jv)w\in K[\text{\rm ad}x_1,\ldots,\text{\rm ad}x_m]$
have no components of degree $\leq1$, we obtain that $[x_i,x_j](h_jt-(g_ju-f_jv)w)=0$ in $L_{m,3}$
and $(\varphi,\phi,\gamma)=1$.

$\text{\rm (iii)}$ Let $m\geq2,c\geq4$ and let $\psi,\phi\in\text{\rm GInn}(L_{m,c})$
be generalized inner automorphisms of the form  
\[
\psi: x_i\to x_i+\sum_{j=1}^{m}[x_i,x_j]f_j, \quad i=1,\ldots,m,
\]
\[
\phi: x_i\to x_i+\sum_{t=1}^{m}[x_i,x_t]g_t, \quad i=1,\ldots,m,
\]
where $f_j,g_t\in K[\text{\rm ad}x_1,\ldots,\text{\rm ad}x_m]$. Then
we know from Lemma \ref{multiplicationginner} that
the composition $\psi\phi$ is of the form
\[
\psi\phi: x_i\to x_i+\sum_{t=1}^{m}[x_i,x_t]g_t+\sum_{j=1}^{m}[x_i,x_j]f_j+
\sum_{j,t=1}^{m}[x_i,x_t,x_j]g_tf_j, \quad i=1,\ldots,m.
\]

Lemma \ref{psi(u)} states that
\[
\psi(u)=u+\sum_{j=1}^{m}[u,x_j]f_j,\quad \phi(u)=u+\sum_{t=1}^{m}[u,x_t]g_t,
\]
for every $u\in L_{m,c}$. Hence
\[
\psi\phi: u\to u+\sum_{t=1}^{m}[u,x_t]g_t+\sum_{j=1}^{m}[u,x_j]f_j+
\sum_{j,t=1}^{m}[u,x_t,x_j]g_tf_j.
\]
If $u$ is an element of the derived algebra $L_{m,c}'$, then
\[
\psi\phi: u\to u+\sum_{t=1}^{m}[u,x_t]g_t+\sum_{j=1}^{m}[u,x_j]f_j+
\sum_{j,t=1}^{m}u(\text{\rm ad}x_t)(\text{\rm ad}x_j)g_tf_j,
\]
where $\text{\rm ad}x_t,\text{\rm ad}x_j,g_t,f_j\in K[\text{\rm ad}x_1,\ldots,\text{\rm ad}x_m]$.
Hence 
\[
\psi\phi(u)=\phi\psi(u), \quad u\in L_{m,c}'.
\]
This means that the commutator
\[
(\psi,\phi)=\psi^{-1}\phi^{-1}\psi\phi\in(\text{\rm GInn}(L_{m,c}),\text{\rm GInn}(L_{m,c}))
\]
of $\psi$ and $\phi$ acts
trivially on $L_{m,c}'$.

Now let us define the generalized normal automorphisms $\rho$ and $\sigma$ in the commutator
subgroup $(\text{\rm GInn}(L_{m,c}),\text{\rm GInn}(L_{m,c}))$
and let $w_1(x_i)=\rho(x_i)-x_i$ and $w_2(x_i)=\sigma(x_i)-x_i$, $i=1,\ldots,m$. Then clearly the elements $w_1(x_i)$
and $w_2(x_i)$ are in  $L_{m,c}'$, i.e. $\rho$ and $\sigma$ act trivially on them. Thus
\begin{align}
&\rho\sigma(x_i)=\rho(x_i+w_2(x_i))=\rho(x_i)+w_2(x_i)=x_i+w_1(x_i)+w_2(x_i)\nonumber\\
&\sigma\rho(x_i)=\sigma(x_i+w_1(x_i))=\sigma(x_i)+w_2(x_i)=x_i+w_1(x_i)+w_2(x_i)\nonumber
\end{align}
which means that $\rho\sigma=\sigma\rho$. Hence $(\text{\rm GInn}(L_{m,c}),\text{\rm GInn}(L_{m,c}))$ is abelian
and so $\text{\rm GInn}(L_{m,c})$ is metabelian.
\end{proof}

\begin{example}
Now we give an explicit proof of the fact that $\text{\rm GInn}(L_{2,3})$ is nilpotent of class 2.
Let $\psi\in\text{\rm GInn}(L_{2,3})$ be a generalized inner automorphism of the form
\begin{align}
\psi(x_1)=&x_1+\alpha[x_1,x_2]+\alpha_1[x_1,x_2,x_1]+\alpha_2[x_1,x_2,x_2]\nonumber\\
\psi(x_2)=&x_2+\beta[x_1,x_2]+\beta_1[x_1,x_2,x_1]+\beta_2[x_1,x_2,x_2]\nonumber
\end{align}
where $\alpha,\alpha_1,\alpha_2,\beta,\beta_1,\beta_2\in K$. Easy calculations give that the inverse
function $\psi^{-1}$ has the form
\begin{align}
\psi^{-1}(x_1)=&x_1-\alpha[x_1,x_2]-(\alpha\beta+\alpha_1)[x_1,x_2,x_1]+(\alpha^2-\alpha_2)[x_1,x_2,x_2]\nonumber\\
\psi^{-1}(x_2)=&x_2-\beta[x_1,x_2]-(\beta^2+\beta_1)[x_1,x_2,x_1]+(\alpha\beta-\beta_2)[x_1,x_2,x_2]\nonumber
\end{align}
If $\phi\in\text{\rm GInn}(L_{2,3})$ is another generalized inner automorphism,
\begin{align}
\phi(x_1)=&x_1+p[x_1,x_2]+p_1[x_1,x_2,x_1]+p_2[x_1,x_2,x_2]\nonumber\\
\phi(x_2)=&x_2+q[x_1,x_2]+q_1[x_1,x_2,x_1]+q_2[x_1,x_2,x_2]\nonumber
\end{align}
where $p,p_1,p_2,q,q_1,q_2\in K$ with inverse
\begin{align}
\phi^{-1}(x_1)=&x_1-p[x_1,x_2]-(pq+p_1)[x_1,x_2,x_1]+(p^2-p_2)[x_1,x_2,x_2]\nonumber\\
\phi^{-1}(x_2)=&x_2-q[x_1,x_2]-(q^2+q_1)[x_1,x_2,x_1]+(pq-q_2)[x_1,x_2,x_2]\nonumber
\end{align}

\noindent calculating the composition $\psi\phi$ we have that
\begin{align}
\psi\phi(x_1)=&x_1+(\alpha+p)[x_1,x_2]+(\alpha_1+p_1-p\beta)[x_1,x_2,x_1]+(\alpha_2+p_2+p\alpha)[x_1,x_2,x_2]\nonumber\\
\psi\phi(x_2)=&x_2+(\beta+q)[x_1,x_2]+(\beta_1+q_1-q\beta)[x_1,x_2,x_1]+(\beta_2+q_2+q\alpha)[x_1,x_2,x_2]\nonumber
\end{align}
and the composition $\phi^{-1}\psi\phi$ is of the form
\begin{align}
\phi^{-1}\psi\phi(x_1)=&x_1+\alpha[x_1,x_2]+(\alpha_1-p\beta+q\alpha)[x_1,x_2,x_1]+\alpha_2[x_1,x_2,x_2]\nonumber\\
\phi^{-1}\psi\phi(x_2)=&x_2+\beta[x_1,x_2]+\beta_1[x_1,x_2,x_1]+(\beta_2-p\beta+q\alpha)[x_1,x_2,x_2].\nonumber
\end{align}
Finally we obtain that $(\psi,\phi)=\psi^{-1}\phi^{-1}\psi\phi$ has the form
\begin{align}
(\psi,\phi)(x_1)=&x_1+(\alpha q-\beta p)[x_1,x_2,x_1]\nonumber\\
(\psi,\phi)(x_2)=&x_2+(\alpha q-\beta p)[x_1,x_2,x_2].\nonumber
\end{align}

Now let $\theta\in\text{\rm GInn}(L_{2,3})$ be a generalized inner automorphism of the form
\begin{align}
\theta(x_1)=&x_1+a[x_1,x_2]+a_1[x_1,x_2,x_1]+a_2[x_1,x_2,x_2]\nonumber\\
\theta(x_2)=&x_2+b[x_1,x_2]+b_1[x_1,x_2,x_1]+b_2[x_1,x_2,x_2]\nonumber
\end{align}
where $a,a_1,a_2,b,b_1,b_2\in K$. Direct calculations give that
\begin{align}
((\psi,\phi),\theta)(x_1)=&x_1+(0.b-0.a)[x_1,x_2,x_1]\nonumber\\
((\psi,\phi),\theta)(x_2)=&x_2+(0.b-0.a)[x_1,x_2,x_2]\nonumber
\end{align}
which means that
\[
((\psi,\phi),\theta)=(\psi,\phi)^{-1}\theta^{-1}(\phi,\psi)\theta=1.
\]
\end{example}

\section{Main Results}
In this section we describe the group of normal automorphisms in terms of
generalized inner automorphisms. We give the explicit
form of the Jacobian matrices of the normal automorphisms and
of the coset representatives of
normally outer $\text{\rm IA}$-automorphisms.

\begin{lemma}\label{c=2} Let $\varphi$ be a normal $\text{\rm IA}$-automorphism of $L_{m,2}$.
Then $\varphi$ is a generalized inner automorphism of $L_{m,2}$. Furthermore $\varphi$ is an inner automorphism
of $L_{m,2}$.
\end{lemma}

\begin{proof} Clearly, $L'_{m,2}$ has a basis $[x_i,x_j]$, $1\leq i<j\leq m$.
Let $\varphi$ be a normal automorphism in $\text{\rm IA}(L_{m,2})$.
If $m=2$, then $\text{\rm IA}(L_{2,2})=\text{\rm Inn}(L_{2,2})$. Since
$\text{\rm Inn}(L_{2,2})\subset \text{\rm N}(L_{2,2})$ then $\varphi$ is an inner automorphism.
In particular $\varphi$ is a generalized inner automorphism.
Let $m\geq 3$ and $\varphi$ be of the form
\begin{align}
\varphi(x_1)=&x_1+[x_1,c_{11}x_1+c_{12}x_2+\cdots+c_{1m}x_m]\nonumber\\
\varphi(x_2)=&x_2+[x_2,c_{21}x_1+c_{22}x_2+\cdots+c_{2m}x_m]\nonumber\\
\vdots & \nonumber\\
\varphi(x_m)=&x_m+[x_m,c_{m1}x_1+c_{m2}x_2+\cdots+c_{mm}x_m]\nonumber,
\end{align}
where $c_{ij}\in K$ for every $i,j=1,\ldots,m$. Now consider the ideal $J_{12}$ of $L_{m,2}$
generated by $x_1+x_2$. As a vector space $J_{12}$ posseses a basis
\[
x_1+x_2,\, [x_1,x_2],\, [x_1+x_2,x_j],\quad j=3,\ldots,m.
\]
Since $\varphi$ is normal $\varphi(x_1+x_2)\in J_{12}$. But
\[
\varphi(x_1+x_2)=x_1+x_2+(c_{12}-c_{21})[x_1,x_2]
+\sum_{j=3}^{m}c_{1j}[x_1+x_2,x_j]+\sum_{j=3}^{m}(c_{2j}-c_{1j})[x_2,x_j],
\]
which means that $\sum_{j=3}^{m}(c_{2j}-c_{1j})[x_2,x_j]\in J_{12}\cap L'_{m,c}$. Then
\[
\sum_{j=3}^{m}d_j[x_2,x_j]=p[x_1,x_2]+\sum_{j=3}^{m}q_j[x_1+x_2,x_j],
\]
for some $p,q_j\in K$, $d_j=c_{2j}-c_{1j}$, $j=3,\dots,m$,
which means that $d_j=0$. Hence $c_{2j}=c_{1j}$, $j=3,\dots,m$.
Similarly, considering the ideals $J_{1k}$ of $L_{m,2}$ generated by $x_1+x_k$ for every $k=3,\ldots,m$, we obtain that
\[
c_{k2}=c_{12},\ldots,c_{k,k-1}=c_{1,k-1},c_{k,k+1}=c_{1,k+1},\ldots,c_{km}=c_{1m}.
\]

Finally, considering the ideals $J_{2k}$ of $L_{m,2}$ generated by $x_2+x_k$, $k=3,\ldots,m$, similar arguments
give that
\[
c_{k1}=c_{21}, \quad k=3,\ldots,m.
\]
Thus 
\[
\varphi=\exp(\text{\rm ad}u),\quad u=c_{21}x_1+c_{12}x_2+c_{13}x_3+\cdots+c_{1m}x_m,
\]
i.e. $\varphi$ is an inner automorphism.
\end{proof}

\begin{lemma}\label{c=3m=2} Every normal $\text{\rm IA}$-automorphism of $L_{2,3}$
is generalized inner.
\end{lemma}

\begin{proof} Let $\varphi$ be a normal $\text{\rm IA}$-automorphism of $L_{2,3}$
such that 
\begin{align}
\varphi(x_1)=&x_1+\alpha[x_1,x_2]+\alpha_1[x_1,x_2,x_1]+\alpha_2[x_1,x_2,x_2]\nonumber\\
\varphi(x_2)=&x_2+\beta[x_2,x_1]+\beta_1[x_2,x_1,x_1]+\beta_2[x_2,x_1,x_2]\nonumber,
\end{align}
where $\alpha,\alpha_1,\alpha_2,\beta,\beta_1\beta_2\in K$.
Let us define $f_1=\beta+\beta_1\text{\rm ad}x_1+\beta_2\text{\rm ad}x_2$ and $f_2=
\alpha+\alpha_1\text{\rm ad}x_1+\alpha_2\text{\rm ad}x_2$. Then we can rewrite $\varphi$ in the following way.
\begin{align}
\varphi(x_1)=&x_1+\sum_{j=1}^{2}[x_1,x_j]f_j,\nonumber\\
\varphi(x_2)=&x_2+\sum_{j=1}^{2}[x_2,x_j]f_j,\nonumber
\end{align}
which completes the proof.
\end{proof}

We know that $\text{\rm Inn}(L_{m,c})\subset\text{\rm N}(L_{m,c})$. If $c=2$, then
$\text{\rm Inn}(L_{m,2})=\text{\rm IN}(L_{m,2})$ by Lemma \ref{c=2}.
But the elements
$\varphi$ of $\text{\rm IN}(L_{m,c})$ are not necessarily inner automorphisms when $c\geq3$.
For example it follows from Lemma \ref{c=3m=2} that
\begin{align}
\varphi(x_1)=&x_1+[x_1,x_2,x_2]\nonumber\\
\varphi(x_2)=&x_2\nonumber
\end{align}
is a normal automorphism which is not an inner automorphism.

\begin{lemma}\label{c>4normal} 
Let $\varphi$ be a normal $\text{\rm IA}$-automorphism of $L_{m,c}$ acting trivially on
$L_{m,c}/L_{m,c}^c$. Then $\varphi$ is a generalized inner
automorphism.
\end{lemma}

\begin{proof}
If $m=2$ then $\varphi$ is of the form
\begin{align}
\varphi(x_1)=&x_1+[x_1,x_2]f_{12},\nonumber\\
\varphi(x_2)=&x_2+[x_2,x_1]f_{21},\nonumber
\end{align}
where $[x_1,x_2]f_{12},[x_2,x_1]f_{21}\in L_{2,c}^c$. This means that
\begin{align}
\varphi(x_1)=&x_1+[x_1,x_1]f_{21}+[x_1,x_2]f_{12},\nonumber\\
\varphi(x_2)=&x_2+[x_2,x_1]f_{21}+[x_2,x_2]f_{12}.\nonumber
\end{align}
Thus $\varphi$ is a generalized inner automorphism.
In the case $c=2,m\geq2$, we know from Lemma \ref{c=2} that such automorphisms are generalized inner.
Hence we can assume that $c\geq 3$,$m\geq 3$.
Let $\varphi$ be a normal automorphism of $L_{m,c}$ acting trivially on
$L_{m,c}/L_{m,c}^c$. Since $c\geq 3$ and $m\geq 3$ then we can assume from Lemma \ref{nia} that
$\varphi$ is of the form
\[
\varphi: x_i\to x_i+\sum_{j=1}^{m}[x_i,x_j]f_{ij},
\]
where $[x_i,x_j]f_{ij}(\text{\rm ad}x_1,\ldots,\text{\rm ad}x_m)$ is in the the center $L_{m,c}^c$ of
the free metabelian nilpotent Lie algebra $L_{m,c}$ for every $i,j=1,\ldots,m$. Such automorphisms form
an abelian subgroup of $\text{\rm Aut}L_{m,c}$.
Let us define the generalized inner automorphism
\[
\varphi_1: x_i\to x_i+\sum_{j=2}^{m}[x_i,x_j]f_{1j}.
\]
Then the composition $\varphi\varphi_1^{-1}$ has the form
\begin{align}
\varphi\varphi_1^{-1}(x_1)=&x_1\nonumber\\
\varphi\varphi_1^{-1}(x_k)=&x_k+[x_k,x_1]f_{k1}+\sum_{j\neq 1,k}^{m}[x_k,x_j](f_{kj}-f_{1j}),\quad k\neq1\nonumber.
\end{align}
Now consider the generalized inner automorphism $\varphi_2: x_i\to x_i+[x_i,x_1]f_{21}$. Then
\begin{align}
\varphi\varphi_1^{-1}\varphi_2^{-1}(x_1)=&x_1\nonumber\\
\varphi\varphi_1^{-1}\varphi_2^{-1}(x_2)=&x_2+\sum_{j=3}^{m}[x_2,x_j]g_{2j}\nonumber\\
\varphi\varphi_1^{-1}\varphi_2^{-1}(x_3)=&x_3+[x_3,x_1]g_{31}+\sum_{j\neq 1,3}[x_3,x_j]g_{3j}\nonumber\\
\vdots&\nonumber\\
\varphi\varphi_1^{-1}\varphi_2^{-1}(x_m)=&x_m+[x_m,x_1]g_{m1}+\sum_{j=2}^{m-1}[x_m,x_j]g_{mj}\nonumber
\end{align}
where $g_{k1}=f_{k1}-f_{21}$ for $k\geq3$ and $g_{kj}=f_{kj}-f_{1j}$ for $k\geq2$, $j\geq2$.
Thus it sufficies to show that $\phi=\varphi\varphi_1^{-1}\varphi_2^{-1}$ is a generalized inner automorphism.
Let $\alpha$ be a nonzero constant and let us consider the ideal $J_{\alpha12}$ of $L_{m,c}$ generated by $\alpha x_1+x_2$.
the vector space $J_{\alpha12}$ modulo $L_{m,c}^3$ has a basis:
\[
\alpha x_1+x_2,\, [x_1,x_2],\, [\alpha x_1+x_2,x_j],\quad j=3,\ldots,m.
\]
Since $\phi$ is normal, $\phi(\alpha x_1+x_2)\in J_{\alpha12}$,
\[
\phi(\alpha x_1+x_2)=\alpha x_1+x_2+\sum_{j=3}^{m}[x_2,x_j]g_{2j},
\]
which means that
\[
\sum_{j=3}^{m}[x_2,x_j]g_{2j}\in J_{\alpha12}\cap L'_{m,c}.
\]
Then
\[
\sum_{j=3}^{m}[x_2,x_j]g_{2j}=[x_1,x_2]P+\sum_{j=3}^{m}[\alpha x_1+x_2,x_j]Q_j,
\]
for some $P,Q_j\in \omega^{c-2}/\omega^{c-1}$, $j=3,\dots,m$. Using the embedding $L_{m,c}$ into the wreath product
we have that
\[
a_2\sum_{j=3}^{m}t_jg_{2j}-\sum_{j=3}^{m}a_jt_2g_{2j}=
a_1(t_2P+\alpha \sum_{j=3}^{m}t_jQ_j)+
a_2(-t_1P+\sum_{j=3}^{m}t_jQ_j)-\sum_{j=3}^{m}a_j(\alpha t_1+t_2)Q_j.
\]
Since $a_1,\dots,a_m$ are free generators of free $K[t_1,\ldots,t_m]/\Omega^{c-1}$-module, for every $j=3,\dots,m$ we have that
\[
t_2g_{2j}=(\alpha t_1+t_2)Q_j.
\]
Thus $\alpha t_1+t_2$ divides $g_{2j}$, $j=3,\dots,m$, for every $\alpha\in K^*$. Since characteristic
of the field $K$ is $0$ we can choose more than $c-2$
distinct scalars $\alpha\in K^*$. Then by nilpotency the function
$g_{2j}(t_1,\ldots,t_m)$ is $0$, $j=3,\dots,m$. Hence
\[
g_{23}=\cdots=g_{2m}=0.
\]
Considering  the ideals $J_{\alpha1k}$, $k=3,\ldots,m$ of $L_{m,c}$ generated by $\alpha x_1+x_k$ the same argument
gives that $g_{kj}=0$, $j\neq1,k$, $k=3,\ldots,m$.

Now let us consider the ideal $J_{\alpha23}$ of $L_{m,c}$ generated by $\alpha x_2+x_3$. 
It has a basis
\[
\alpha x_2+x_3,\, [x_2,x_3],\, [\alpha x_2+x_3,x_j],\quad j\neq2,3,
\]
modulo $L_{m,c}^3$. Since $\phi$ is normal, $\phi(\alpha x_2+x_3)\in J_{\alpha23}$.
\[
\phi(\alpha x_2+x_3)=\alpha x_2+x_3+\sum_{j\neq3}[x_3,x_j]g_{3j}.
\]
This means that $[x_3,x_1]g_{31}\in J_{\alpha23}\cap L'_{m,c}$ because we know that
$g_{3k}=0$, $k\neq1,3$. Then
\[
[x_3,x_1]g_{31}=[x_2,x_3]P+\sum_{j\neq2,3}[\alpha x_2+x_3,x_j]Q_j,
\]
for some $P,Q_j\in \omega^{c-2}/\omega^{c-1}$, $j\neq2,3$. Using the embedding $L_{m,c}$ into the wreath product,
considering only the coefficient of $a_1$ we have that
\[
t_3g_{31}=(\alpha x_2+x_3)Q_1.
\]
Thus $\alpha t_2+t_3$ divides $g_{31}$ for every $\alpha\in K^*$. Hence the function
$g_{31}(t_1,\ldots,t_m)$ is $0$ and $g_{31}=0$.

Finally, considering  the ideals $J_{\alpha 2k}$, $k=4,\ldots,m$ of $L_{m,c}$ generated by $\alpha x_2+x_k$ the same argument
gives that $g_{k1}=0$, $k=4,\ldots,m$. Hence $\phi=\varphi\varphi_1^{-1}\varphi_2^{-1}=1$, i.e.
$\varphi=\varphi_2\varphi_1$ which means that $\varphi$ is a generalized inner automorphism.
\end{proof}

\begin{theorem}\label{IAnormal} 
Let $\varphi$ be a normal $\text{\rm IA}$-automorphism of $L_{m,c}$.
Then $\varphi$ is a generalized inner automorphism.
\end{theorem}

\begin{proof}
We argue by induction on the nilpotency class $c$ of $L_{m,c}$. If $c=2$, 
the result follows from Lemma \ref{c=2} (in this case, each normal $\text{\rm IA}$-automorphism
is inner). Now consider a normal $\text{\rm IA}$-automorphism $\varphi$ of $L_{m,c}$, $c>2$, of the form
\[
\varphi: x_i\to x_i+\sum_{j=1}^{m}[x_i,x_j](f_{ij,0}+\cdots+f_{ij,c-2}), 
\]
where $f_{ij,0}\in K$, $f_{ij,k}\in\omega^k/\omega^{k+1}$, $k=1,\ldots,c-2$.
Then $\varphi$ induces a normal $\text{\rm IA}$-automorphism on $L_{m,c}/L_{m,c}^c$. By induction, since this
quotient is isomorphic to $L_{m,c-1}$, there exists a generalized inner automorphism $\psi:L_{m,c}\to L_{m,c}$ such that 
\[
\varphi(x_i)=\psi(x_i)+\sum_{j=1}^{m}[x_i,x_j]f_{ij,c-2},\quad i=1,\ldots,m.
\]
It follows for $i=1,\ldots,m$ that 
\begin{align}
\psi^{-1}\varphi(x_i)&=x_i+\sum_{j=1}^{m}\psi^{-1}([x_i,x_j]f_{ij,c-2})\nonumber\\
&=x_i+\sum_{j=1}^{m}[x_i,x_j]f_{ij,c-2}.\nonumber
\end{align}
Thus $\phi=\psi^{-1}\varphi$ is a normal $\text{\rm IA}$-automorphism
of $L_{m,c}$ acting trivially on $L_{m,c}/L_{m,c}^c$. By Lemma \ref{c>4normal}, $\phi$ is a generalized inner automorphism,
and so is $\varphi=\psi\phi$.
\end{proof}

Now we give one of the main results which is obtained as a direct consequence of Lemma \ref{normalLmc},
Lemma \ref{c=2}, Lemma \ref{c=3m=2} and Theorem \ref{IAnormal}.

\begin{corollary}\label{structure} Let $K^*$ denote the set of invertible elements of the field $K$. Then
\begin{align}
\text{\rm (i)}&\, \text{\rm N}(L_{m,1})\cong K^*;\nonumber\\
\text{\rm (ii)}&\, \text{\rm N}(L_{2,2})\cong K^*\rightthreetimes\text{\rm Inn}(L_{2,2});\nonumber\\
\text{\rm (iii)}&\, \text{\rm N}(L_{2,3})\cong K^*\rightthreetimes\text{\rm GInn}(L_{2,3});\nonumber\\
\text{\rm (iv)}&\, \text{\rm N}(L_{m,c})={\rm GInn}(L_{m,c}), \quad m\geq3,c\geq2\, \text{\rm or}\, m=2,c\geq4,\nonumber
\end{align}
where $\rightthreetimes$ stands for the semi-direct product of the groups.
\end{corollary}

Now we describe the group structure of the group of normal automorphisms $\text{\rm N}(L_{m,c})$.

\begin{theorem}\label{nautgroup}
$\text{\rm (i)}$ The group $\text{\rm N}(L_{m,2})$, $m\geq3$, is abelian;

$\text{\rm (ii)}$ The group $\text{\rm N}(L_{m,3})$, $m\geq3$, is nilpotent of class 2;

$\text{\rm (iii)}$ The group $\text{\rm N}(L_{m,c})\in\mathfrak{A}^2$, $m\geq2,c\geq4$ or $(m,c)=(2,2)$, is metabelian;

$\text{\rm (iv)}$ The group $\text{\rm N}(L_{2,3})\in\mathfrak{N}_2\mathfrak{A}$, is nilpotent of class two--by--abelian.
\end{theorem}

\begin{proof} Let $\psi_{\alpha},\phi_{\beta}\in\text{\rm N}(L_{2,2})\cong K^*\rightthreetimes\text{\rm Inn}(L_{2,2})$
be normal automorphisms of the form
\begin{align}
\psi_{\alpha}(x_1)=&\alpha x_1+\alpha \alpha_2[x_1,x_2]\nonumber\\
\psi_{\alpha}(x_2)=&\alpha x_2+\alpha \alpha_1[x_2,x_1]\nonumber
\end{align}
\begin{align}
\psi_{\beta}(x_1)=&\beta x_1+\beta \beta_2[x_1,x_2]\nonumber\\
\psi_{\beta}(x_2)=&\beta x_2+\beta \beta_1[x_2,x_1]\nonumber
\end{align}
where $\alpha_1,\alpha_2,\beta_1,\beta_2\in K$ and $\alpha,\beta\in K^*$. Easy calculations give that
\begin{align}
\psi_{\alpha}^{-1}(x_1)=&\alpha^{-1} x_1-\alpha^{-2} \alpha_2[x_1,x_2]\nonumber\\
\psi_{\alpha}^{-1}(x_2)=&\alpha^{-1} x_2-\alpha^{-2} \alpha_1[x_2,x_1]\nonumber
\end{align}
\begin{align}
\psi_{\beta}^{-1}(x_1)=&\beta^{-1} x_1-\beta^{-2} \beta_2[x_1,x_2]\nonumber\\
\psi_{\beta}^{-1}(x_2)=&\beta^{-1} x_2-\beta^{-2} \beta_1[x_2,x_1]\nonumber
\end{align}

By direct calculations we obtain that the commutator
$(\psi_{\alpha},\phi_{\beta})=\psi_{\alpha}^{-1}\phi_{\beta}^{-1}\psi_{\alpha}\phi_{\beta}$ has the form
\[
(\psi_{\alpha},\phi_{\beta}): x_i\to x_i+\sum_{j=1}^{m}
(\alpha^{-1}\alpha_j(\beta^{-1}-1)+\beta^{-1}\beta_j(1-\alpha^{-1}))[x_i,x_j], \quad i=1,2.
\]
This means that $(\psi_{\alpha},\phi_{\beta})\in\text{\rm GInn}(L_{2,2})=\text{\rm Inn}(L_{2,2})$
which is abelian from Theorem \ref{gigroup}. Hence
$\text{\rm N}(L_{2,2})$
is metabelian.

We know from Corollary \ref{structure} that if
$m\geq3$ or $m=2,c\geq4$ then $\text{\rm N}(L_{m,c})={\rm GInn}(L_{m,c})$. Applying Theorem \ref{gigroup}
we get that $\text{\rm N}(L_{m,2})$ is abelian when  $m\geq3$,
$\text{\rm N}(L_{m,3})$ is nilpotent of class 2 when $m\geq3$
and that $\text{\rm N}(L_{m,c})$ is metabelian when  $m\geq2,c\geq4$.
Thus it remains to show that
$\text{\rm N}(L_{2,3})$ is a nilpotent of class two--by--abelian group.

Now let $\psi_{\alpha},\phi_{\beta}$ in $\text{\rm N}(L_{2,3})$ be normal automorphisms of the form
\[
\psi_{\alpha}: x_i\to \alpha x_i+f_i, \quad i=1,\ldots,m,
\]
\[
\phi_{\beta}: x_i\to \beta x_i+g_i, \quad i=1,\ldots,m,
\]
where $f_i,g_i\in L_{m,c}'$ and $\alpha,\beta\in K^*$. Clearly the inverse functions are of the form
\[
\psi_{\alpha}^{-1}: x_i\to \alpha^{-1} x_i+f'_i, \quad i=1,\ldots,m,
\]
\[
\phi_{\beta}^{-1}: x_i\to \beta^{-1} x_i+g'_i, \quad i=1,\ldots,m,
\]
where $f'_i,g'_i\in L_{m,c}'$. Easy calculations give that the commutator $(\psi_{\alpha},\phi_{\beta})$ of $\psi_{\alpha}$ and $\phi_{\beta}$
is included in $\text{\rm GInn}(L_{2,3})$ which is
nilpotent of class $2$ by Theorem \ref{gigroup}. Thus $\text{\rm N}(L_{2,3})$ is a nilpotent of class two--by--abelian group.
\end{proof}

Now we have collected necessary information for the description of the group of normally outer
automorphisms $\Gamma\text{\rm N}(L_{m,c})$.
We shall find the coset representatives of the normal subgroup $\text{\rm IN}(L_{m,c})$
of the group $\text{\rm IA}(L_{m,c})$ of IA-automorphisms $L_{m,c}$, i.e.,
we shall find a set of IA-automorphisms $\theta$ of $L_{m,c}$ such that
the factor group $\text{\rm I}\Gamma\text{\rm N}(L_{m,c})=\text{\rm IA}(L_{m,c})/\text{\rm IN}(L_{m,c})$
of the outer IA-automorphisms of $L_{m,c}$
is presented as the disjoint union of the cosets
$\text{\rm IN}(L_{m,c})\theta$.

\begin{lemma}\label{m=2}
Let $m=2$, then the group of normally outer $\text{\rm IA}$-automorphisms $\text{\rm I}\Gamma\text{\rm N}(L_{2,c})$ is trivial.
\end{lemma}

\begin{proof}
Let $\varphi$ be an $\text{\rm IA}$-automorphism of $L_{2,c}$.
Then $\varphi$ has the form
\begin{align}
\varphi:&x_1\to x_1+[x_1,x_2]f\nonumber\\
&x_2\to x_2+[x_1,x_2]g,\nonumber
\end{align}
where $f,g \in K[\text{\rm ad}x_1,\text{\rm ad}x_2]$. Then clearly
\begin{align}
\varphi:&x_1\to x_1+[x_1,x_1]f_1+[x_1,x_2]f_2\nonumber\\
&x_2\to x_2+[x_2,x_1]f_1+[x_2,x_2]f_2,\nonumber
\end{align}
where $f_1=g,f_2=f$, i.e. $\varphi$ is a generalized inner automorphism or from Theorem \ref{IAnormal}
$\varphi$ is a normal $\text{\rm IA}$-automorphism. Thus $\text{\rm IA}(L_{2,c})=\text{\rm IN}(L_{2,c})$.
\end{proof}

\begin{theorem}\label{nonnormalJ}
\text{\rm (i)} Let $\varphi$ be a normal $\text{\rm IA}$-  automorphism of the form
\[
\varphi: x_i\to x_i+\sum_{j=1}^{m}[x_i,x_j]f_j, \quad i=1,\ldots,m,
\]
where $f_j\in K[\text{\rm ad}x_1,\ldots,\text{\rm ad}x_m]$.
Then the Jacobian matrix of $\varphi$ is
\[
J(\varphi)=I_m+\left(\begin{array}{cccc}
t_2f_2+\cdots+t_mf_m&-t_2f_1&\cdots&-t_mf_1\\
-t_1f_2&\sum_{j\neq2}t_jf_j&\cdots&-t_mf_2\\
-t_1f_3&-t_2f_3&\cdots&-t_mf_3\\
 \vdots&\vdots&\ddots&\vdots\\
-t_1f_m&-t_2f_m&\cdots&\sum_{j\neq m}t_jf_j\\
\end{array}\right),
\]

\text{\rm (ii)} Let $\Theta$ be the set of automorphisms $\theta$ of $L_{m,c}$ with Jacobian matrix of the form
\[
J(\theta)=I_m+\left(\begin{array}{cccc}
0&f_{12}(\hat t_2)&\cdots&f_{1m}\\
p_2(\hat t_1)&f_{22}&\cdots&f_{2m}\\
p_3(\hat t_1)&f_{32}&\cdots&f_{3m}\\
\vdots&\vdots&\ddots&\vdots\\
p_m(\hat t_1)&f_{m2}&\cdots&f_{mm}\\
\end{array}\right),
\]
where $p_i,f_{ij}$, are polynomials of degree $\leq c-1$ without constant terms with
the following conditions
\[
\sum_{i=2}^mt_ip_i\equiv 0,\quad \sum_{i=1}^mt_if_{ij}\equiv 0
\quad (\text{\rm mod }\Omega^{c+1}), \quad j=2,\ldots,m,
\]
$p_i=p_i(\hat t_1)$, $i=1,\ldots,m$, does not depend on $t_1$,
and  $f_{12}=f_{12}(\hat t_2)$ does not depend on $t_2$.

Then $\Theta$ consists of coset representatives of the subgroup $\text{\rm IN}(L_{m,c})$ 
of the group $\text{\rm IA}(L_{m,c})$ and $\text{\rm I}\Gamma\text{\rm N}(L_{m,c})$ is a disjoint union of
the cosets $\text{\rm IN}(L_{m,c})\theta$, $\theta\in \Theta$.

\text{\rm (iii)} Let $\Psi$ be the set of normal $\text{\rm IA}$-automorphisms $\psi$ of $L_{m,c}$ with Jacobian matrix of the form
\[
J(\psi)=I_m+\left(\begin{array}{cccc}
\sum_{j\neq 1}t_jq_j(T_j)&-t_2q_1(T_1)&\cdots&-t_mq_1(T_1)\\
-t_1q_2(T_2)&\sum_{j\neq2}t_jq_j(T_j)&\cdots&-t_mq_2(T_2)\\
-t_1q_3(T_3)&-t_2q_3(T_3)&\cdots&-t_mq_3(T_3)\\
\vdots&\vdots&\ddots&\vdots\\
-t_1q_m(T_m)&-t_2q_m(T_m)&\cdots&\sum_{j\neq m}t_jq_j(T_j)\\
\end{array}\right),
\]
where $q_j(T_j)$, $j=1,\ldots,m$, are polynomials of degree $\leq c-1$ in $\Omega^2$ with
the following conditions
\[
\sum_{i=2}^mq_j(T_j)\equiv 0\quad (\text{\rm mod }\Omega^{c+1}),
\]
and $q_j(T_j)$ depends on $t_j,\ldots,t_m$ only, $j=1,\ldots,m$.

Then $\Psi$ consists of coset representatives of the subgroup $\text{\rm Inn}(L_{m,c})$ 
of the group $\text{\rm IN}(L_{m,c})$ and $\text{\rm IN}(L_{m,c})/\text{\rm Inn}(L_{m,c})$ is a disjoint union of
the cosets $\text{\rm Inn}(L_{m,c})\psi$, $\psi\in \Psi$.
\end{theorem}

\begin{proof} $\text{\rm (i)}$ 
Let $\varphi$ be a normal $\text{\rm IA}$-automorphism of the form
\[
\psi: x_i\to x_i+\sum_{j=1}^{m}[x_i,x_j]f_j, \quad i=1,\ldots,m,
\]
where $f_j\in K[\text{\rm ad}x_1,\ldots,\text{\rm ad}x_m]$.
The Jacobian matrix of $\varphi$ is
\[
J(\varphi)=\left(\frac {\partial \varphi({x_j})}{\partial x_i}\right)
=\left(\begin{matrix}
\frac {\partial\varphi({x_1})}{\partial x_1}&\cdots&\frac {\partial \varphi({x_m})}{\partial x_1}\\
\vdots&\ddots&\vdots\\
\frac {\partial\varphi({x_1})}{\partial x_m}&\cdots&\frac {\partial \varphi({x_m})}{\partial x_m}\\
\end{matrix}\right)\in M_m(K[t_1,\ldots,t_m]/\Omega^c).
\]
Easy calculations give

\[
\frac{\partial\varphi(x_j)}{\partial x_i}=\delta_{ij}+\begin{cases}
\sum_{r\neq j}t_rf_r&i=j,\\
-t_jf_i&i\neq j,\end{cases}
\]
where $\delta_{ij}$ is Kronecker symbol. Thus we obtain the desired form of the matrix $J(\varphi)$.

$\text{\rm (ii)}$
When $m=2$ then from Lemma \ref{m=2} the factor group
$\text{\rm IA}(L_{2,c})/\text{\rm IN}(L_{2,c})$ is trivial which satisfies the conditions. Let $m\geq 3$.
Since $\text{\rm Inn}(L_{m,c})$ is included in
the group of normal automorphisms, the factor group $\text{\rm IA}(L_{m,c})/\text{\rm IN}(L_{m,c})$ is the homomorphic image
of  $\text{\rm IA}(L_{m,c})/\text{\rm Inn}(L_{m,c})$. Then from Lemma \ref{formula of out}
we can consider the Jacobian matrix of the $\text{\rm IA}$-automorphism
$\psi$ of the form
\[
J(\psi)=I_m+\left(\begin{array}{llll}
s(t_2,\ldots,t_m)&f_{12}&\cdots&f_{1m}\\
t_1q_2(t_2,t_3,\ldots,t_m)+r_2(t_2,\ldots,t_m)&f_{22}&\cdots&f_{2m}\\
t_1q_3(t_3,\ldots,t_m)+r_3(t_2,\ldots,t_m)&f_{32}&\cdots&f_{3m}\\
\ \ \ \ \ \ \ \vdots&\ \ \vdots&\ \ddots&\ \ \vdots\\
t_1q_m(t_m)+r_m(t_2,\ldots,t_m)&f_{m2}&\cdots&f_{mm}\\
\end{array}\right),
\]
where $s,q_i,r_i,f_{ij}$ are polynomials of degree $\leq c-1$ without constant terms with the conditions
\[
s+\sum_{i=2}^mt_iq_i\equiv 0,\quad \sum_{i=2}^mt_ir_i\equiv 0,\quad \sum_{i=1}^mt_if_{ij}\equiv 0
\quad (\text{\rm mod }\Omega^{c+1}), \quad j=2,\ldots,m,
\]
$s=s(t_2,\ldots,t_m)$, $r_i=r_i(t_2,\ldots,t_m)$, $i=1,\ldots,m$, does not depend on $t_1$, $q_i(t_i,\ldots,t_m)$,
$i=2,\ldots,m$, does not depend on $t_1,\ldots,t_{i-1}$
and  $f_{12}$ does not contain a summand $dt_2$, $d\in K$.

Let
\[
f_1=0,\quad
f_k=q_k,\quad k=2,\ldots,m,
\]
and let us define the normal automorphism
\[
\varphi: x_i\to x_i+\sum_{j=1}^{m}[x_i,x_j]f_j, \quad i=1,\ldots,m.
\]
Then from $\text{\rm (i)}$ the Jacobian matrix of $\varphi$ is of the form
\[
J(\varphi)=I_m+\left(\begin{array}{cccc}
-s&0&\cdots&0\\
-t_1q_2&-s-t_2q_2&\cdots&-t_mq_2\\
-t_1q_3&-t_2q_3&\cdots&-t_mq_3\\
 \vdots&\vdots&\ddots&\vdots\\
-t_1q_m&-t_2q_m&\cdots&-s-t_mq_m\\
\end{array}\right).
\]
Let us denote the $m\times 2$ matrix consisting of the first two columns of $J(\varphi\psi)$ and $I_m$ by
$J(\varphi\psi)_2$ and $I_{m2}$, respectively. Direct calculations give that $J(\varphi\psi)_2$ is of the form
\[
J(\varphi\psi)_2=I_{m2}+\left(\begin{array}{cccc}
-s^2&-sf_{12}+f_{12}\\
-s(t_1q_2+r_2)+r_2&\ast\\
-s(t_1q_3+r_3)+r_3&\ast\\
\vdots&\vdots\\
-s(t_1q_m+r_m)+r_m&\ast\\
\end{array}\right),
\]
where we have denoted by $\ast$ the corresponding entries of the second column of
the Jacobian matrix of $\varphi\psi$.

Now let
\[
g_1=0,\quad
g_k=-sq_k,\quad k=2,\ldots,m,
\]
and let us define the normal automorphism
\[
\phi: x_i\to x_i+\sum_{j=1}^{m}[x_i,x_j]g_j, \quad i=1,\ldots,m.
\]
The Jacobian matrix of $\phi$ is of the form
\[
J(\phi)=I_m+\left(\begin{array}{cccc}
s^2&0&\cdots&0\\
st_1q_2&s(s+t_2q_2)&\cdots&st_mq_2\\
st_1q_3&st_2q_3&\cdots&st_mq_3\\
 \vdots&\vdots&\ddots&\vdots\\
st_1q_m&st_2q_m&\cdots&s(s+t_mq_m)\\
\end{array}\right).
\]
Calculating $J(\phi\varphi\psi)$ we have that
\[
J(\phi\varphi\psi)_2=I_{m2}+\left(\begin{array}{cccc}
-s^4&f_{12}(1-s+s^2-s^3)\\
-s^3t_1q_2+r_2(1-s+s^2-s^3)&\ast\\
-s^3t_1q_3+r_3(1-s+s^2-s^3)&\ast\\
\vdots&\vdots\\
-s^3t_1q_m+r_m(1-s+s^2-s^3)&\ast\\
\end{array}\right).
\]
Repeating this process sufficiently many times, we get that the $(1,1)$-th entry and the coefficients of
the elements $t_1q_j$, $j=2,\ldots,m$, are zero, because $L_{m,c}$ is nilpotent. So we have the form
\[
J(\gamma)_2=I_{m2}+\left(\begin{array}{cccc}
0&g_{12}\\
p_2(\hat t_1)&\ast\\
p_3(\hat t_1)&\ast\\
\vdots&\vdots\\
p_m(\hat t_1)&\ast\\
\end{array}\right),
\]
where $p_i=p_i(\hat t_1)$, $i=2,\ldots,m$, does not depend on $t_1$, $g_{12}$ does not contain a summand $dt_2$, $d\in K$.
Let us express $g_{12}$ as
\[
g_{12}=t_2f+\hat f_2,
\]
where $\hat f_2$ does not depend on $t_2$ and $f\in\Omega$ because $g_{12}$ does not contain a summand $dt_2$, $d\in K$. Let us consider the normal automorphism
\[
\phi_1: x_i\to x_i+[x_i,x_1]f, \quad i=1,\ldots,m.
\]
The Jacobian matrix of $\phi_1$ is of the form
\[
J(\phi_1)=I_m+\left(\begin{array}{cccc}
0&-t_2f&\cdots&-t_mf\\
0&t_1f&\cdots&0\\
0&0&\cdots&0\\
\vdots&\vdots&\ddots&\vdots\\
0&0&\cdots&t_1f\\
\end{array}\right),
\]

Calculating $J(\phi_1\gamma)$ we have that
\[
J(\phi_1\gamma)_2=I_{m2}+\left(\begin{array}{cccc}
0&t_1f(t_2f+\hat f_2)+t_2f+\hat f_2-t_2f\\
t_1fp_2+p_2&\ast\\
t_1fp_3+p_3&\ast\\
\vdots&\vdots\\
t_1fp_m+p_m&\ast\\
\end{array}\right).
\]

Now let
\[
g_1=0,\quad
g_k=fp_k,\quad k=2,\ldots,m,
\]
and let us define the normal automorphism
\[
\phi_2: x_i\to x_i+\sum_{j=1}^{m}[x_i,x_j]g_j, \quad i=1,\ldots,m.
\]

Calculating $J(\phi_2\phi_1\gamma)$ we see that the summands $-t_1fp_j$ in the first column disappears:
\[
J(\phi_2\phi_1\gamma)_2=I_{m2}+\left(\begin{array}{cccc}
0&t_1t_2f^2+t_1f\hat f_2+\hat f_2\\
p_2&\ast\\
p_3&\ast\\
\vdots&\vdots\\
p_m&\ast\\
\end{array}\right).
\]
Let us consider the $(1,2)$-th entry $t_1t_2f^2+t_1f\hat f_2+\hat f_2$ of the matrix $J(\phi_2\phi_1\gamma)$ and express the element $f$ as
\[
f=t_2F+\hat F_2,
\] 
where $\hat F_2$ does not depend on $t_2$. Now we have that
\begin{align}
t_1t_2f^2+t_1f\hat f_2+\hat f_2&=t_1t_2(f^2+F\hat f_2)+(t_1\hat F_2+1)\hat f_2\nonumber\\
&=t_1t_2h+\hat h_2,\nonumber
\end{align}
where $\hat h_2=(t_1\hat F_2+1)\hat f_2$ does not depend on
$t_2$ and $h=f^2+F\hat f_2$.
Note that the minimal degree of the monomials of the summand which depend on $t_2$ (in this step this is $t_1t_2h$)
is bigger than of the minimal degree of the corresponding summand $t_2f$ of the previous step which means that the degree increases.

We repeat the process one more step and consider the normal automorphism
\[
\varphi_1: x_i\to x_i+[x_i,x_1](\text{\rm ad}x_1h), \quad i=1,\ldots,m.
\]
Calculating $J(\varphi_1\phi_2\phi_1\gamma)$ we have that
\[
J(\varphi_1\phi_2\phi_1\gamma)_2=I_{m2}+\left(\begin{array}{cccc}
0&t_1^2h(t_1t_2h+\hat h_2)+t_1t_2h+\hat h_2-t_2t_1h\\
t_1^2hp_2+p_2&\ast\\
t_1^2hp_3+p_3&\ast\\
\vdots&\vdots\\
t_1^2hp_m+p_m&\ast\\
\end{array}\right).
\]

Now let
\[
g_1=0,\quad
g_k=t_1hp_k,\quad k=2,\ldots,m,
\]
and let us define the normal automorphism
\[
\varphi_2: x_i\to x_i+\sum_{j=1}^{m}[x_i,x_j]g_j, \quad i=1,\ldots,m.
\]
Then
\[
J(\varphi_2\varphi_1\phi_2\phi_1\gamma)_2=I_{m2}+\left(\begin{array}{cccc}
0&t_1^3t_2h^2+t_1^2h\hat h_2+\hat h_2\\
p_2&\ast\\
p_3&\ast\\
\vdots&\vdots\\
p_m&\ast\\
\end{array}\right).
\]

Let us consider the $(1,2)$-th entry $t_1^3t_2h^2+t_1^2h\hat h_2+\hat h_2$ of the matrix $J(\varphi_2\varphi_1\phi_2\phi_1\gamma)_2$ and express the element $h$ as
\[
h=t_2H+\hat H_2,
\] 
where $\hat H_2$ does not depend on $t_2$. Now we have that
\begin{align}
t_1^3t_2h^2+t_1^2h\hat h_2+\hat h_2 &=t_1^2t_2(t_1h^2+H\hat h_2)+(t_1^2\hat H_2+1)\hat h_2 \nonumber\\
&=t_1^2t_2Q(t_1,\ldots,t_m)+\hat Q(\hat t_2),\nonumber
\end{align}
Again, the length of the summands in $t_1^2t_2Q(t_1,\ldots,t_m)$ which depend on $t_2$
increases step by step. Repeating this argument sufficiently many times, by nilpotency, we get finally that
\[
J(\gamma)_2=I_{m2}+\left(\begin{array}{cccc}
0&q(\hat t_2)\\
p_2(\hat t_1)&\ast\\
p_3(\hat t_1)&\ast\\
\vdots&\vdots\\
p_m(\hat t_1)&\ast\\
\end{array}\right).
\]

Hence, starting from an arbitrary coset of IA-automorphisms $\text{\rm IN}(L_{m,c})\psi$,
we have found that it contains an automorphism $\theta\in\Theta$ with Jacobian matrix
prescribed in the theorem.
Now, let $\theta_1$ and $\theta_2$ be two different automorphisms in $\Theta$ with
$\text{\rm IN}(L_{m,c})\theta_1=\text{\rm IN}(L_{m,c})\theta_2$. Hence,
there exists a nontrivial automorphism $\varphi$ in $\text{\rm IN}(L_{m,c})$
such that $\theta_1=\varphi\theta_2$. Direct calculations show that this
is in contradiction with the form of $J(\theta_1)$.

$\text{\rm (iii)}$
Let $\varphi$ be a normal $\text{\rm IA}$-automorphism of $L_{m,c}$.
From $\text{\rm (i)}$, the Jacobian matrix of $\varphi$ is
\[
J(\varphi)=I_m+\left(\begin{array}{cccc}
t_2f_2+\cdots+t_mf_m&-t_2f_1&\cdots&-t_mf_1\\
-t_1f_2&\sum_{j\neq2}t_jf_j&\cdots&-t_mf_2\\
-t_1f_3&-t_2f_3&\cdots&-t_mf_3\\
 \vdots&\vdots&\ddots&\vdots\\
-t_1f_m&-t_2f_m&\cdots&\sum_{j\neq m}t_jf_j\\
\end{array}\right),
\]
where $f_j(t_1,\ldots,t_m)\in K[t_1,\ldots,t_m]$, $j=1,\ldots,m$. When $c=2$ then from Lemma \ref{c=2}, $\text{\rm IN}(L_{m,2})=\text{\rm Inn}(L_{m,2})$. As a result we may consider that $f_j(t_1,\ldots,t_m)\in \Omega$.

Let us express the polynomials $f_j(t_1,\ldots,t_m)$, $j=2,\ldots,m$, in the following way:
\[
f_j(t_1,\ldots,t_m)=t_1\overline f_j(t_1,\ldots,t_m)+h_j(T_2).
\]
Now let $u\in L_{m,c}$ be of the form
\[
u=-\sum_{i>1}[x_i,x_1]\overline f_j(\text{\rm ad}x_1,\ldots,\text{\rm ad}x_m), \quad g_{i1}(t_1,\ldots,t_m)\in \Omega,
\]
and let consider the inner automorphism $\phi_1=\exp(\text{\rm ad}u)$.
Then the Jacobian matrix of $\phi_1$ has the form
\[
J(\phi_1)=I_m+\left(
\begin{array}{cccc}
-t_1G_1&-t_2G_1&\cdots&-t_mG_1\\
-t_1 G_2&-t_2G_2&\cdots&-t_mG_2\\
\vdots&\vdots&\ddots&\vdots\\
-t_1G_m&-t_2G_m&\cdots&-t_mG_m
\end{array}
\right),
\]
where
\begin{align}
G_1=&t_2\overline f_2-t_3\overline f_3-\cdots-t_m\overline f_m\nonumber\\
G_2=&-t_1\overline f_2,\, G_3=-t_1\overline f_3,\ldots,\,G_m=-t_1\overline f_m.\nonumber
\end{align}

The element $u$ belongs to the commutator ideal of $L_{m,c}$ and
the linear operator $\text{\rm ad}u$ acts trivially on $L'_{m,c}$.
Hence $\exp(\text{\rm ad}u)$ is the identity map restricted on $L'_{m,c}$.
Since the automorphism $\varphi$ is IA, we obtain that
\[
J(\phi_1\varphi)_2=I_{m2}+\left(\begin{array}{cccc}
t_2h_2(T_2),+\cdots+t_mh_m(T_2),&-t_2F_1\\
-t_1h_2(T_2),&t_1F_1+\sum_{j=3}^mt_jh_j(T_2),\\
-t_1h_3(T_2),&-t_2h_3(T_2)\\
 \vdots&\vdots\\
-t_1h_m(T_2),&-t_2h_m(T_2)\\
\end{array}\right),
\]

Now we write $h_i(T_2)$ in the form
\[
h_i(T_2)=t_2h_i'(T_2)+h_i''(T_3),\quad i=3,\ldots,m,
\]
and define
\[
\phi_2=\exp(\text{\rm ad}u_2),\quad
u_2=\sum_{i=3}^m[x_i,x_2]h_i'(\text{\rm ad}x_2,\ldots,\text{\rm ad}x_m).
\]
Then we obtain that
\[
J(\phi_2\phi_1\phi_0\psi)_2=\left(\begin{array}{cccc}
1+t_2H_2(T_2)+\cdots+t_mh_m''(T_m)&-t_2F'_1\\
-t_1H_2(T_2)&\ast\\
-t_1h_3''(T_3)&\ast\\
\vdots&\vdots\\
-t_1h_m''(T_3)&\ast\\
\end{array}\right),
\]
\[
H_2(T_2)=h_2(T_2)-\sum_{i=3}^mt_ih_i'(T_2).
\]
Repeating this process we construct inner automorphisms $\phi_3,\ldots,\phi_{m-1}$ such that
\[
\psi=\phi_{m-1}\cdots\phi_2\phi_1\varphi,
\]
\[
J(\phi_{m-1}\cdots\phi_2\phi_1\varphi)_2=\left(\begin{array}{cccc}
1+t_2H_2(T_2)+\cdots+t_mH_m(T_m)&-t_2H_1(T_1)\\
-t_1H_2(T_2)&\ast\\
-t_1H_3(T_3)&\ast\\
\vdots&\vdots\\
-t_1H_m(T_m)&\ast\\
\end{array}\right),
\]

Hence, starting from an arbitrary coset of normal $\text{\rm IA}$-automorphisms $\text{\rm Inn}(L_{m,c})\varphi$,
we found that it contains an automorphism $\psi\in\Psi$ with Jacobian matrix
prescribed in the theorem.
Now, let $\psi_1$ and $\psi_2$ be two different automorphisms in $\Psi$ with
$\text{\rm Inn}(L_{m,c})\psi_1=\text{\rm Inn}(L_{m,c})\psi_2$. Hence,
there exists a nonzero element $u\in L_{m,c}$
such that $\psi_1=\exp(\text{\rm ad}u)\psi_2$. Direct calculations show that this
is in contradiction with the form of $J(\psi_1)$.
\end{proof}

\begin{example}
When $m=3$ the results of Theorem \ref{nonnormalJ}
have the following simple form. 
If $\varphi$ is a normal automorphism of the form
\begin{align}
\varphi:&x_1\to x_1+[x_1,x_2]f_2+[x_1,x_3]f_3\nonumber\\
&x_2\to x_2+[x_2,x_1]f_1+[x_2,x_3]f_3\nonumber\\
&x_3\to x_3+[x_3,x_1]f_1+[x_3,x_2]f_2\nonumber
\end{align}
where $f_1,f_2,f_3\in K[\text{\rm ad}x_1,\text{\rm ad}x_2,\text{\rm ad}x_3]$
then the Jacobian matrix of $\varphi$ is
\[
J(\varphi)=\left(\begin{array}{cccc}
1+t_2f_2+t_3f_3&-t_2f_1&-t_3f_1\\
-t_1f_2&1+t_1f_1+t_3f_3&-t_3f_2\\
-t_1f_3&-t_2f_3&1+t_1f_1+t_2f_2\\
\end{array}\right).
\]

The Jacobian matrix of the normally outer automorphism $\theta$ is
\[
J(\theta)=\left(\begin{array}{cccc}
1&f_{12}(t_1,t_3)&f_{13}\\
t_3p(t_2,t_3)&1+f_{22}&f_{23}\\
-t_2p(t_2,t_3)&f_{32}&1+f_{33}\\
\end{array}\right),
\]
where $p(t_2,t_3),f_{ij}$, are polynomials of degree $\leq c-1$ without constant terms with
the following conditions
\[
t_1f_{1j}+t_2f_{2j}+t_3f_{3j}\equiv 0
\quad (\text{\rm mod }\Omega^{c+1}), \quad j=2,3,
\]
$p(t_2,t_3)$ does not depend on $t_1$
and  $f_{12}=f_{12}(t_1,t_3)$ does depend on $t_2$.
\end{example}

\section*{Acknowledgements}

The author is very thankful to Professor Vesselin Drensky for many useful suggestions and discussions.
He is very thankful to Professor Ualbai Umirbaev for the useful discussions on the Freiheitssatz
for Lie and associative algebras.
Finally, he is grateful to the Institute of Mathematics and Informatics of
the Bulgarian Academy of Sciences for the creative atmosphere and the warm hospitality
and during his visit when
most of this project was carried out.

\end{document}